\documentclass[preprint,12pt]{elsarticle}
\usepackage{amsfonts}
\usepackage{mathrsfs}

\usepackage{amsmath}
\usepackage{cases}
\usepackage{hyperref}
\usepackage{amssymb}
\usepackage{graphicx}
\usepackage{subfigure}
\usepackage{slashbox}
\usepackage{pict2e}

\usepackage{amsthm}

\biboptions{numbers,sort&compress}

\begin{document}
\begin{frontmatter}

\title{On the  method of directly defining inverse mapping \\ for nonlinear differential equations}


\author{Shijun Liao $^{a,b,c}$ \corref{cor1} }

\ead{sjliao@sjtu.edu.cn} \cortext[cor1]{Corresponding author.}

\author{Yinlong Zhao $^a$}

\address{

$^a$ State Key Laboratory of Ocean Engineering\\
School of Naval Architecture, Ocean and Civil Engineering\\
Shanghai Jiao Tong University, Shanghai 200240, P.R. China\\

\vspace{0.3cm}

$^b$ Collaborative Innovation Center for Advanced Ship and Deep-Sea Exploration (CISSE), Shanghai 200240,  P.R. China

\vspace{0.3cm}

$^c$ MOE Key Laboratory in Scientific and Engineering Computing\\
Dept. of Mathematics, Shanghai Jiao Tong University, Shanghai 200240, P.R. China

}

\begin{abstract}
In scientific computing,  it is  time-consuming to calculate an inverse operator ${\mathscr A}^{-1}$ of a differential equation ${\mathscr A}\varphi = f$, especially when ${\mathscr A}$ is a highly nonlinear operator.   In this paper,  based on the homotopy analysis method (HAM),   a new approach, namely the method of directly defining inverse mapping (MDDiM), is proposed to gain analytic approximations of nonlinear differential equations.   In other words, one can solve a nonlinear differential equation  ${\mathscr A}\varphi = f$  by means of directly defining an inverse mapping $\mathscr J$, i.e. without calculating any inverse operators.      Here, the inverse mapping $\mathscr J$ is even unnecessary to  be  explicitly  expressed in a differential form, since  ``mapping'' is a more general  concept than ``differential operator''.    To guide how to directly define an inverse mapping $\mathscr J$,  some rules are provided.   Besides,  a convergence theorem is proved, which guarantees that a convergent series solution given by the MDDiM must be a solution of problems under consideration.    In addition, three nonlinear differential equations are used  to illustrate the validity and potential of the MDDiM, and especially the great freedom and large flexibility of directly defining inverse mappings for various types of nonlinear problems.  The method of directly defining inverse mapping (MDDiM) might open a completely new, more general way to solve nonlinear problems in science and engineering, which is fundamentally different from traditional methods.
\end{abstract}

\begin{keyword}
Homotopy analysis method \sep  analytical approximation    \sep nonlinear differential equation \sep direct definition of inverse mapping
\end{keyword}

\end{frontmatter}


\section{Motivation}

For a differential equation ${\mathscr A}\varphi = f$, where ${\mathscr A}$ is a differential operator and $f$ is a known function,  one can quickly gain its common solution $u={\mathscr A}^{-1} f$, when the inverse operator ${\mathscr A}^{-1}$ is known,or if  it  is  easy to gain.  Unfortunately,  lots of CPU times (i.e. a large amount of money) are often consumed to calculate inverse operator ${\mathscr A}^{-1}$ in general.

Can we solve nonlinear differential equations by means of directly defining an inverse mapping, i.e. without calculating any inverse operators?  If so,  lots of CPU time (and money) can be saved.  This is the motivation of this work.

Traditionally,  perturbation techniques \cite{Nayfeh} are widely used to gain analytic approximations of a nonlinear differential equation ${\mathscr A} \varphi =  f$.  If there exists a small physical parameter $\epsilon$, and besides  if the nonlinear operator $\mathscr A$ contains a linear ones, i.e. $\mathscr{A} = \mathscr{L} + \mathscr{N}$, one can express \[ \varphi = \varphi_0 + \varphi_1 \epsilon + \varphi_2 \epsilon^2+\cdots \] and transfer the original nonlinear equation ${\mathscr A} \varphi =  f$  into  an infinite number of linear sub-problms
\[    {\mathscr L} [\varphi_0] = f, \;  {\mathscr L}[\varphi_m] = Q_m(\varphi_0,\varphi_1,\cdots,\varphi_{m-1}), \;\; m=1,2,3,\cdots,\]
where $Q_m$ is dependent upon the known terms $\varphi_0,\varphi_1,\cdots,\varphi_{m-1}$ and thus is known.   Note that these linear sub-problems have a {\em close}  relationship with the original equation: they use the {\em same}  linear  operator $\mathscr{L}$ that is the linear part of the original equation ${\mathscr A}\varphi = f$.    In applied mathematics, there exist  many  methods  that transfer a nonlinear problem into a series of linear sub-problems.  Traditionally,  these  linear sub-problems often have rather close relationship with the original ones, but are often  difficult to  solve, because it is generally time-consuming  to obtain an inverse operator even for a  linear equation.  Sometimes,  the linear part even does {\em not} contain the highest order of derivatives so that the linear sub-problems become ``singular'' since there are more boundary/initial conditions.   It is even worse, when ${\mathscr A}$  does {\em not} contain any linear parts at all!   This is  mainly because  perturbation techniques and other traditional analytic approximation methods can {\em not} provide  us {\em freedom}  to  choose the related linear operators of linear sub-problems, that determine their equation-types.

Fortunately, we have such kind of freedom in the frame of the Homotopy Analysis Method (HAM)  \cite{liaoPhd, liaobook1, liaobook2, liaobook3}, an analytic approximation technique for highly nonlinear problems.   Based on homotopy, a basic concept in topology about continuous variation,  the HAM can easily transfer a nonlinear equation into an infinite number of linear sub-problems.   Compared to perturbation techniques and other traditional analytic methods, the HAM has some advantages.   First of all,  the HAM can transfer a nonlinear problem into an infinite number of linear sub-problems {\em without} any small/large physical parameters.   In other words, the HAM works even if there  do {\em not} exist any small/large physical parameters in governing equations and initial/boundary conditions!   Secondly, the HAM provides us great freedom to choose the equation-type of the linear sub-problems, i.e. the freedom to choose an auxiliary linear operator ${\mathscr L}$ for linear sub-problems, even if the original nonlinear operator $\mathscr A$ does {\em not} contain any linear parts, since we have great freedom in the frame of the HAM to construct different homotopies (or variations).   Especially,  unlike perturbation techniques and other analytic methods,  the HAM provides us a convenient way to guarantee the convergence of solution series by means of introducing the so-called ``convergence-control parameter'' into the solution series.  With these advantages, the HAM has been widely applied to solve nonlinear problems in lots of fields  \cite{liaobook1, liaobook2, liaobook3, Abbasbandy2007PLA, Tao2007JCR, Yongyan2009WaveMotion, KV2009CNSNS, Mastroberardino2011-CNSNS, KV2012, Motsa2010B,  Liang2014}.  For example,    the HAM was successfully applied to give, for the first time,  the theoretical prediction of the so-called steady-state resonant waves (with time-independent spectrum) in deep and finite depth of water \cite{xu2012JFM, Liu2014JFM} for full wave equations, which was confirmed in 2015 by the physical experiments \cite{Liu2015JFM}.  For details, please refer to \cite{BookChap3-Liao2015}.   This illustrates the potential and novelty of the HAM, since a truly new method should bring us something new and different!

Here,  it should be emphasized that the HAM provides us great freedom to choose the equation-type and auxiliary linear operator $\mathscr L$ of the linear sub-problems.  Such kind of freedom is so large that,  in the frame of the HAM, a 2nd-order Gelfand equation can be  solved very easily  by means of  transferring it into an infinite number of 4th-order (two-dimensional)  or 6th-order (three-dimensional) linear differential equations, and the convergent series solutions were in  good agreement to numerical ones,  as illustrated by Liao and Tan \cite{Liao2007}.   Note that it was traditionally believed that a 2nd-order differential equation could be replaced {\em only} by an infinite number of linear differential equations at the {\em same} order, if  perturbation techniques \cite{Nayfeh},  Adomian Decomposition Method \cite{Adomian1994} and other traditional methods are used.  So, Liao and Tan's approach \cite{Liao2007} is difficult to understand from the traditional viewpoints, although it works quite well in practice.   However, this simple example in \cite{Liao2007} reveals something novel and unusual of the HAM:  the HAM can provide us freedom to directly define the auxiliary linear operator $\mathscr L$ of linear sub-problems.   Obviously, if we can transfer a nonlinear equation into an infinite number of linear sub-problems whose inverse linear operators are {\em known} or {\em directly defined}, it becomes straight-forward  to  solve the original nonlinear problem.   This is indeed true:  in this paper, we generalize the HAM-based analytic approach in \cite{Liao2007}  and propose the ``method of directly defining inverse operator'' (MDDiM) in the frame of the HAM.

Can we solve a nonlinear differential equation by means of {\em directly defining an inverse operator}?  This is an open question up to now, to the best of our knowledge.    A  positive  answer  is  given  in  the frame of the HAM \cite{liaoPhd, liaobook1, liaobook2, liaobook3} in this paper.  The paper is organized as follows.  In \S~2, the method of directly defining inverse mapping (MDDiM) and a theorem of convergence are briefly described.  In \S~3, we give the detailed derivation of the MDDiM and prove the theorem of convergence mentioned in \S~2.  In \S~4, three examples are used to illustrate how to apply the MDDiM to solve nonlinear differential equations.   Some discussions and concluding remarks are given in \S~5.

\section{Method of directly defining inverse mapping (MDDiM)}

First of all,   we briefly describe the basic ideas of the method of directly defining inverse mapping (MDDiM).

Let us consider a $n$th-order nonlinear differential equation
\begin{equation}
  \mathscr{N}[u(x)]=0,\;\;\; x\in\Omega,   \label{geq:original}
\end{equation}
subject to the $   \mu  $ linear boundary conditions
 \begin{equation}
 \mathscr{B}_i [u] = \beta_i, \;\;\;\; \mbox{at $x=\alpha_i$}, \;\;\;\; i =1 ,2,3, \cdots,    \mu  ,   \label{bc:original}
 \end{equation}
where $u(x)$ is a unknown function, $x$ is an independent-variable, $\Omega$ is an interval of $x$,  $\mathscr{N}$ denotes a nonlinear operator,  ${\mathscr B}_i$ is a linear operator,  $1 \leq     \mu   \leq n$ are positive integers, $\alpha_i \in \Omega$ and $\beta_i$ ($1\leq i\leq \mu$) are constants,  respectively.  Note that $n=\mu$ for linear problems, but this is unnecessary for nonlinear ones.

Let
\[   {S}_\infty  =  \left\{  \varphi_1(x), \varphi_2(x), \cdots   \right\} \]
denote a complete set of  an  infinite  number of  base functions that are
\emph{linearly independent}.
All functions that are expressed by $S_\infty$ form a set of functions, denoted by
\begin{equation}
  V = \left\{ \left. \sum_{k=1}^{+\infty} a_k \varphi_k(x)  \right| a_k \in \mathbb{R}\right\}.  \label{def:set:V}
\end{equation}
   Besides, let
\[   S^* =  \left\{  \varphi_1(x), \varphi_2(x), \cdots, \varphi_\mu   \right\}  \]
denote a set, consist of the first $\mu$ simplest base functions of $S_\infty$.
All functions that are expressed by $S^*$ form a set of functions, denoted by
\begin{equation}
 V^* = \left\{ \left. \sum_{k=1}^{\mu} a_k \varphi_k(x)  \right| a_k \in \mathbb{R}\right\}. \label{def:set:V:star}
\end{equation}
Assume that  $u(x) \in V$ and the $\mu$ unknown coefficients of the expression
\[  u^*(x) = \sum_{n=1}^{\mu} a_n \; \varphi(x)  \in V^*  \]
can be uniquely determined by the $\mu$ linear boundary conditions (\ref{bc:original}), i.e.
\[   {\mathscr B}_i \left[\sum_{n=1}^{\mu} a_n \; \varphi(x) \right] = \beta_i, \;\;\;\; \mbox{at $x=\alpha_i$}, \;\;\;\; i =1 ,2,3, \cdots,    \mu.   \]
Then, we call $u^* \in V^*$ the primary solution.
Write
\begin{equation}
 \hat{S} = \left\{  \varphi_{\mu+1}(x), \varphi_{\mu+2}(x), \cdots   \right\}. \label{def:S:hat:general}
 \end{equation}
All functions that are expressed by $\hat{S}$ form a set of functions, denoted by
\begin{equation}
  \hat{V} = \left\{ \left. \sum_{k=\mu+1}^{+\infty} b_k \varphi_k(x)  \right| b_k \in \mathbb{R}\right\}. \label{def:set:V:hat}
\end{equation}
Obviously,  $V=\hat{V} \cup V^*$.
Similarly,  let
\begin{equation}\label{V_R}
S_R  =  \{\psi_1(x),\psi_2(x),\cdots \}
\end{equation}
 be an infinite set of base functions that are
linearly independent,  and all functions expressed by $S_R$ form  a  set  of functions, denoted by
\begin{equation}
U = \left\{\sum_{m=1}^{+\infty} c_m \psi_m(x) \left| c_m \in \mathbb{R} \right. \right\}.  \label{def:set:U}
\end{equation}
Assume that   $ {\mathscr N}[u(x)] \in U$, say, the nonlinear differential operator  ${\mathscr N}$ is a kind of mapping from $V$ to $U$, i.e. ${\mathscr N}: V \to U$.

In the frame of the MDDiM, the series solution of $u(x)$ is given by
\begin{equation}
u(x) = u_0(x) + \sum_{k=1}^{+\infty} u_k(x) ,   \label{def:series-solution}
\end{equation}
where $u_0(x)$ is an {\em initial guess} that satisfies all linear boundary conditions (\ref{bc:original}), and besides we have great freedom to choose it. Here,  $u_k(x)\in V$ is given by
\begin{equation}
u_k(x) = \chi_k \; u_{k-1}(x) + c_0 \; {\mathscr J}\left[ \delta_{k-1}(x)\right]+ \sum^{\mu}_{n=1}a_{k,n} \varphi_n(x),  \label{def:delta-mapping}
\end{equation}
with the definitions
\begin{eqnarray}
\delta_n(x) &=& \mathscr{D}_n \left\{ {\mathscr N}\left[u_0(x)+\sum_{j=1}^{+\infty} u_j(x) \; q^j\right] \right\}  \in U,  \label{def:delta[n]}
\end{eqnarray}
and
\begin{equation} \label{def:chi}
\chi_{k}=\left\{
\begin{aligned}
         0, \;\; & k\leq 1, \\
         1, \;\; & k>1,
\end{aligned} \right.
\end{equation}
where $c_0$ is the so-called ``convergence-control parameter'', which we have great freedom to choose,  $\mathscr J$ is a directly defined inverse mapping,  the operator $\mathscr{D}_n$ is the so-called $n$th-order homotopy-derivative, defined by
\begin{equation}
\mathscr{D}_n \phi  = \frac{1}{n!}\left.  \frac{\partial^n  \phi}{\partial q^n} \right|_{q=0},  \label{def:D[n]}
\end{equation}
whose properties were proved by Liao \cite{Liao2009} and are briefly listed in the Appendix of this paper.
Note that we can regard
\[   \hat{u}_k(x) = \chi_k \; u_{k-1}(x) + c_0 \; {\mathscr J}\left[ \delta_{k-1}(x)\right]  \]
as a special solution of $u_k(x)$, and
\[  u_k^*(x) = \sum^{\mu}_{i=1}a_{k,i} \varphi_i(x) \in V^* \]
 as a primary solution of $u_k(x)$, respectively.   According to (\ref{def:D[n]}),   $\delta_n$ defined by (\ref{def:delta[n]})  can be regarded as the coefficient of Maclaurin series of the governing equation with respect to the embedding parameter $q\in[0,1]$,  say,
\begin{equation}
{\mathscr N}\left[ \sum_{n=0}^{+\infty} u_n (x) q^n \right] = \sum_{n=0}^{+\infty}\delta_n(x) \; q^n.
\end{equation}
This provides us a simple way to gain $\delta_n(x)$  for $n\geq 0$.

  In (\ref{def:delta-mapping}),  the unknown coefficients $a_i$ $(1\leq i \leq  \mu)$  are determined by the $\mu$ boundary conditions
\begin{equation}
{\mathscr B}_i\left[ u_k(x)-\chi_k \; u_{k-1}(x)  \right] = c_{i}\; \Delta_{i,k-1}(x), \;\; \mbox{at $x=\alpha_i$, $1\leq i\leq \mu$, }   \label{bc:mth}
\end{equation}
with the definition
\begin{eqnarray}
\Delta_{i,n}(x) &=& {\mathscr B}_i [u_n(x)] -(1-\chi_{n+1})\beta_i,
\end{eqnarray}
where $c_i$ ($1\leq i\leq \mu$) is the so-called ``convergence-control parameters'', which we have great freedom to choose.

In (\ref{def:delta-mapping}),   $\mathscr J$ denotes a directly {\em defined} mapping from $U\to V$, with the following rules:
\begin{enumerate}
  \item[(I)] $\mathscr{J}$ is linear, i.e
  \begin{equation*}
   \forall \alpha,\beta \in \mathbb{R}, \forall x, y\in U, \mathscr{J}(\alpha x+\beta y)=\alpha \mathscr{J}(x) + \beta \mathscr{J}(y);
  \end{equation*}
  \item[(II)] $\mathscr{J}$ is  injective, say, the kernel of $\mathscr{J}$ is $\{\mathbf{0}\}$, i.e
  \begin{equation*}
    \{x | x \in U, \mathscr{J}(x) = \mathbf{0}\}=\{\mathbf{0}\};
  \end{equation*}
  \item[(III)]  $\mathscr{J}[\delta_{m}(x)]$ contains each base function $\varphi_i\in \hat{S}$ $( \mu+1\leq i < +\infty)$  as $m\to +\infty$ ;
  \item[(IV)] $\mathscr{J}$ is finite, i.e.  there exists  such a finite constant $K$ that for any $\varphi\in V$  it holds  \[  \frac{||{\mathscr J}[{\mathscr N}[\varphi]]||}{||\varphi||}  \leq  K.    \]
\end{enumerate}
Here,  $V, \hat{S}, U $ are defined by (\ref{def:set:V}), (\ref{def:S:hat:general}) and (\ref{def:set:U}), respectively.

It should be emphasized that there exist an auxiliary parameter $c_0$ in (\ref{def:delta-mapping}) and the $\mu$ auxiliary parameters $c_1, c_2, \cdots, c_\mu$ in the boundary conditions (\ref{bc:mth}).   All of them have no physical meanings, but in theory we have great freedom to choose their values.   Mostly, if $c_0$ and $c_i$ ($1\leq i\leq \mu$) are properly chosen,  we can guarantee the convergence of the series solution (\ref{def:series-solution}),  as illustrated later.   This is the reason why we call  $c_0, c_1, c_2, \cdots, c_\mu$   ``the convergence-control parameters''.

In addition,  the following theorem can be proved.

{\bf Theorem of Convergence}.  {\em If the convergence-control parameters $c_0$, $c_1$, $\cdots$, $c_{\mu}$ and the directly defined  inverse mapping  ${\mathscr J}$ are properly chosen so that the series (\ref{def:series-solution}) is absolutely convergent, then it must be a solution of the original equation (\ref{geq:original}) and (\ref{bc:original}). }

According to the above theorem,   we only need choose proper mapping $\mathscr J$ and proper convergence-control parameters $c_0, c_1, c_2, \cdots, c_\mu$ so as to guarantee the convergence of solution series.   The  proof of this convergence theorem and the detailed derivation of the  MDDiM will be given below.

\section{The detailed derivations of MDDiM}

The above-mentioned MDDiM is based on the homotopy analysis method (HAM) \cite{liaoPhd,liaobook1,liaobook2,liaobook3}, a analytic approximation technique for highly nonlinear differential equations.

The HAM is based on homotopy, a fundamental concept in topology, which describes a continuous variation (or deformation) between an initial guess and an exact solution of an equation.   Without loss of generality, let us take the nonlinear differential equation (\ref{geq:original}) and (\ref{bc:original}) as an example.   Let $u_0(x)$ denote an initial guess of the solution $u(x)$ that satisfies the $\mu$ linear boundary conditions (\ref{bc:original}), $c_0, c_1, c_2, \cdots, c_\mu$ are the $(\mu+1)$ non-zero auxiliary parameters  without physical meanings (called ``convergence-control parameters''),  $\mathscr{L}:V\to U$ is an auxiliary linear operator with the property ${\mathscr L}[0]=0$,   and  $q\in[0,1]$ an embedding parameter of homotopy, respectively.  To build a continuous  variation  (or deformation), denoted by $\Phi(x;q)$,  from the initial guess $u_0(x)$ to the exact solution $u(x)$, we construct the so-called zeroth-order deformation equation
\begin{equation}\label{def:zeroth deformation}
    (1-q)\mathscr{L}[\Phi(x;q)-u_0(x)]=q \;
    c_0  \; \mathscr{N}[\Phi(x;q)],   \;\;\; q \in [0,1],
\end{equation}
subject to the $\mu$ linear boundary conditions
 \begin{equation}
 (1-q) {\mathscr B}_i \left[ \Phi(x;q) -u_0(x) \right] =  q \; c_i \;\left\{  {\mathscr B}_i[\Phi(x;q)] -  \beta_i \right\}, \;\;\;\; \mbox{at $x=\alpha_i$},   \label{bc:zero}
 \end{equation}
where  $1\leq i\leq \mu$.    Obviously,  $\Phi(x;0) = u_{0}(x)$ when $q=0$, since ${\mathscr L}[0]=0$.
Besides,  $\Phi(x;1) = u(x)$ when $q=1$, since $c_i\neq 0$ for $0\leq i \leq \mu$.    In other words,  Eqs. (\ref{def:zeroth deformation}) and (\ref{bc:zero})  define  a continuous variation  $\Phi(x;q)$  from the initial guess  $ u_{0}(x)$ to the solution $u(x)$ of the original equations  (\ref{geq:original}) and (\ref{bc:original}), as the homotopy  parameter $q$ increases from $0$ to $1$.
Assuming that the solution $\Phi(x;q)$ is \emph{analytic} at $q=0$, the Maclaurin series of $\Phi(x;q)$ with respect to $q$ reads
\begin{equation}\label{def:series:q}
    \Phi(x;q)=u_{0}(x)+\sum^{+\infty}_{k=1}u_{k}(x)q^k,
\end{equation}
where
\begin{equation}\label{def:uk}
  u_k(x)=\left.\frac{1}{k!}\frac{\partial^k \Phi(x;q)}{\partial q^k}\right|_{q=0} = \mathscr{D}_k(\Phi).
\end{equation}
Here, $\mathscr{D}_k$ is called the $k$th-order  homotopy-derivative operator, defined by (\ref{def:D[n]}).   For  properties  and  theorems  about  $\mathscr{D}_k$ in details,  please refer to \cite{Liao2009} and \S~4.2 of Liao's book \cite{liaobook2}.

Applying the $k$th-order  homotopy-derivative operator  $\mathscr{D}_k$
to both sides of the zeroth-order deformation equations  (\ref{def:zeroth deformation})  and (\ref{bc:zero}),   it is straightforward to obtain the $k$th-order deformation equation
\begin{equation}\label{def:mth}
    \mathscr{L}[u_{k}(x)-\chi_k u_{k-1}(x)]=c_{0} \;    \delta_{k-1}(x), \; k\geq 1
\end{equation}
subject to the $   \mu  $ linear boundary conditions
\begin{equation}
{\mathscr B}_i\left[ u_k(x) -\chi_k u_{k-1}(x)\right] = c_i  \; \Delta_{i,k-1}(x), \;\;\;\; \mbox{at $x=\alpha_i$}, \; 1\leq i \leq \mu  ,   \label{bc:mth}
\end{equation}
where $\chi_n$ is defined by (\ref{def:chi}), and
\begin{eqnarray}
  \delta_{n}(x) &=& \mathscr{D}_{n}\{\mathscr{N}[\Phi(x;q)]\}, \label{def:Rm} \\
  \Delta_{i,n}(x) &=& \mathscr{D}_{n}\left\{  {\mathscr B}_i[\Phi(x;q)] -  \beta_i \right\}
  = {\mathscr B}_i[u_n(x)] -  (1-\chi_{n+1})\beta_i.
\end{eqnarray}
Note that $\delta_{k-1}(x)$ and $\Delta_{i,n}(x)$ are only dependent upon $u_{0}(x)$, $\cdots$, $u_{k-1}(x)$ and thus are {\em known} for the unknown term $u_k(x)$.   So,  $u_{k}(x)$ is determined by the {\em linear} differential equation (\ref{def:mth}) with the $\mu$ {\em linear} boundary conditions (\ref{bc:mth}).

It should be emphasized here that, in the frame of the HAM,  one has great {\em freedom} to choose the
auxiliary linear operator $\mathscr{L}$, the initial guess
$u_{0}(x)$,  and especially the so-called convergence-control parameters $c_{0}$ and $c_1, c_2, \cdots, c_\mu$, as pointed out by Liao \cite{Liao2007}.  Assuming that
all of them are {\em properly} chosen so that the Maclaurin series
(\ref{def:series:q}) converges at $q=1$, one gets the series solution
\begin{equation}\label{def:series:q=1}
    u(x)=u_{0}(x)+\sum^{+\infty}_{k=1}u_{k}(x).
\end{equation}
The $m$th-order approximation of $u(x)$ reads
\begin{equation}\label{def:Um}
 u (x) \approx u_0(x) +  \sum^{m}_{k=1}u_{k}(x).
\end{equation}
Thus, in essence, the HAM transfers a nonlinear problem into an infinite number of linear sub-problems.  However, unlike perturbation methods \cite{Nayfeh}, we do {\em not} need any small/large {\em physical} parameters at all in the frame of the HAM for such kind of transformation.   In addition, unlike perturbation methods \cite{Nayfeh}, we have now great freedom to choose the auxiliary linear operator $\mathscr L$.  More importantly,  the so-called ``convergence-control parameters'' $c_0$ and $c_1, c_2, \cdots, c_\mu$ provide a convenient  way to guarantee the convergence of the solution series, as illustrated by lots of  successful applications of the HAM \cite{liaobook1, liaobook2, liaobook3,KV2012}.

 \subsection{Normal strategy of the HAM}

In the frame of the HAM,  normally,  one often {\em chooses} such a proper auxiliary linear operator $\mathscr L$ that the linear high-order deformation equations (\ref{def:mth}) and (\ref{bc:mth}) are {\em easy} to solve, and  besides that the convergence of the solution series is guaranted  by means of choosing proper convergence-control parameters $c_0$ and $c_1, c_2, \cdots, c_\mu$.   This is mainly because we have great {\em freedom} to choose $\mathscr L$ and the convergence-control parameters  in the frame of the HAM.   This is completely different from perturbation techniques.   To guide how to choose $\mathscr L$,  Liao \cite{liaobook1, liaobook2} suggested a few rules described below.

Assume that
\begin{equation}
u(x) = \sum_{m=1}^{+\infty} b_m \; \varphi_m(x) \in V,   \label{solution-expression}
\end{equation}
where $V$ is defined by (\ref{def:set:V}).   We call it ``the solution expression'' of $u(x)$, which plays an important role in the normal frame of the HAM.  Unlike perturbation methods, the solution expression is the starting point of the HAM, since it greatly influences the choice of the auxiliary linear operator $\mathscr L$.   As suggested by Liao \cite{liaobook1, liaobook2},   $\mathscr L$ should be chosen in such a way that
\begin{enumerate}
\item[(a)] there exists a unique solution $u_k(x)$ of the $k$th-order deformation equation  ({\em Rule of Solution Existence});

\item[(b)]  $u_k(x) \in V$ ({\em Rule of Solution Expression});

\item[(c)] $\sum\limits_{k=0}^{+\infty} u_k(x)$ contains all base functions.  ({\em Rule of Completeness}).
\end{enumerate}

In addition, due to the Rule of Solution Existence,  $\mathscr L$ should be chosen in such a way that it holds
\begin{equation}
{\mathscr L} [\varphi_i(x)] = 0, \;\;\; \forall  \varphi_i(x) \in V^*,  \;\;\; 1\leq i\leq \mu,
\end{equation}
and
\begin{equation}
 {\mathscr L}[\varphi_i(x)] \neq 0, \;\;\; \forall \varphi_i(x) \in \hat{V},\;\;\; i >\mu,
\end{equation}
where $V^*$ and $\hat{V}$  are defined by (\ref{def:set:V:star}) and (\ref{def:set:V:hat}), respectively, since there exist the $\mu$ linear boundary conditions (\ref{bc:original}).   In other words,  $\varphi_i \in V^*$, here $1\leq i\leq \mu$, is a primary solution of ${\mathscr L}[u(x)] = 0$.   Let ${\mathscr L}^{-1}: U\to \hat{V}$ denote the inverse operator of $\mathscr L$, where $U$ and $\hat{V}$ are defined by (\ref{def:set:U}) and (\ref{def:set:V:hat}), respectively.  We have the common solution
\begin{equation}\label{u_k:solution:L}
u_{k}(x)=\chi_k u_{k-1}(x)+c_0 {\mathscr L}^{-1} [\delta_{k-1}(x)] + \sum^{\mu}_{i=1}a_{k,i} \varphi_i(x)
\end{equation}
of the high-order deformation equation (\ref{def:mth}), where the unknown coefficients $a_{k,i}$ ($1\leq i\leq \mu$) are uniquely determined by the $\mu$ linear boundary conditions (\ref{bc:mth}).

In essence, the key of this normal strategy of the HAM is to gain the inverse operator ${\mathscr L}^{-1}$ of the auxiliary linear operator $\mathscr L$.  Unfortunately,  it is often {\em time-consuming} to gain an inverse operator ${\mathscr L}^{-1}$ of a differential equation, unless the linear operator $\mathscr L$ is simple enough.   Due to this restriction,  we often had to  choose very simple auxiliary linear operators ${\mathscr L}$ in the frame of the HAM.    This  widely  restricts  applications of the HAM.  To overcome this limitations,  a new strategy of the HAM  is suggested below.

\subsection{New strategy of the HAM}

Write  $\mathscr{J} = {\mathscr L}^{-1}:U\to\hat{V}$, which is an inverse linear operator of $\mathscr L$.   It should be emphasized that,  in the frame of the HAM, we have great {\em freedom} to choose $\mathscr L$.  In theory, it means that we have great {\em freedom} to {\em directly} choose  ${\mathscr L}^{-1}$, i.e.  we also have great freedom to {\em define} ${\mathscr J}:U \to \hat{V}$, {\em directly}, without choosing the auxiliary linear operator $\mathscr L$ at all!

  Then, the solution $u_{k}(x)$ of Eq. (\ref{def:mth}) reads
\begin{equation}\label{u_k:solution:L}
  u_{k}(x)=\chi_k u_{k-1}(x)+c_0 \mathscr{J}[\delta_{k-1}(x)] + \sum^{\mu}_{i=1}a_{k,i} \varphi_i(x),
\end{equation}
where
\[  \hat{u}_k(x) = \chi_k u_{k-1}(x)+c_0 \mathscr{J}[\delta_{k-1}(x)] \]
is a special solution of $u_k(x)$,
\[ u_k^*(x) = \sum^{\mu}_{i=1}a_{k,i} \varphi_i(x) \]
is a primary solution of $u_k(x)$,  and $a_{k,1}, a_{k,2},\dots,a_{k,\mu}$ are constants to be uniquely determined by the $\mu$ linear boundary condition (\ref{bc:mth}), respectively.    Here, it should be emphasized that, according to  (\ref{u_k:solution:L}), it is {\em unnecessary} to know the specific form of the auxiliary linear operator $\mathscr{L}: V\to U$,  since  the inverse operator $\mathscr{J}:U \to \hat{V}$ is {\em defined}  directly.   In this way, it is unnecessary to spend any CPU times to calculate the inverse operator $\mathscr{J}$, since it is know!

The new strategy of the HAM is fundamentally different from the normal ones.
In the normal HAM,  one should first choose (or define) a proper (but simple enough) auxiliary linear operator $\mathscr L$,  then solve the linear high-order deformation equation (\ref{def:mth}), say,  find out its inverse operator ${\mathscr J} = {\mathscr L}^{-1}$ by means of spending lots of CPU times.   This is often {\em time-consuming} and sometimes even {\em impossible}, especially when $\mathscr L$ is complicated.    However,  using the new strategy of the HAM,  one can  neglect the auxiliary linear operator $\mathscr L$ completely, but {\em define} the inverse linear operator ${\mathscr J} =  {\mathscr L}^{-1}$  {\em directly}!  In this way,  the high-order deformation equation can be quickly solved, since it is {\em unnecessary} to gain the inverse operator ${\mathscr L}^{-1}$ at all!

It should be emphasized that it is the HAM that provides us great freedom to choose the auxiliary linear operator $\mathscr L$, so that we further have the great freedom to {\em directly} {\em define} its inverse operator ${\mathscr J}={\mathscr L}^{-1}$.   For simplicity, we call this approach  ``the method of directly defining inverse mapping'' (MDDiM).

\subsection{Some rules of directly defining the  inverse mapping $\mathscr J$}

Like the normal strategy of the HAM, the initial guess $u_0(x)$, the primary solutions and the inverse operator $\mathscr J$ should be chosen in such a way that
\begin{enumerate}
\item[(A)] there exists a unique solution $u_k(x)$ of the $k$th-order deformation equation  ({\em Rule of Solution Existence});

\item[(B)]  $u_k(x) \in V$ ({\em Rule of Solution Expression});

\item[(C)] $\sum\limits_{k=0}^{+\infty }u_k(x)$ contains all base functions({\em Rule of Completeness}).
\end{enumerate}

First of all, to obey the ``Rule of Solution Expression'',  we should choose an initial guess $u_0(x) \in V$.  Since  we have great freedom to choose $u_0(x)$ in the frame of the HAM, we can choose
\[  u_0(x) =\sum_{i=1}^{\mu} a_{0,i} \; \varphi_i(x) \in V^*, \]
where $V^*$ is defined by (\ref{def:set:V:star}), and the coefficients $a_{0,i}$ ($1\leq i\leq \mu$) are determined by the linear boundary conditions (\ref{bc:original}).

Secondly,  since the linear differential equation (\ref{def:mth}) has the $\mu$ linear boundary conditions (\ref{bc:mth}),   the new strategy should provide the $\mu$ primary solutions of it.  Obviously, to obey ``the Rule of Solution Expression'', each primary solution $u^*_{k}(x)$ must belong to $V^*$.  Thus,  we {\em directly define} the primary solution
\[   {u}_k^*(x) = \sum_{i=1}^{\mu} {a}_{k,i} \; \varphi_i(x) \in V^*,    \]
where  $a_{k,i}$  are  the  unknown  constants, which can be  determined by the linear boundary conditions (\ref{bc:mth}).

Thirdly, to obey the ``Rule of Solution Expression'',  we should have
\[  \delta_{k-1}(x)\in U,  {\mathscr J}[\delta_{k-1}(x)] \in \hat{V},\]
for $k\geq 1$,  and  the special solution $\hat{u}_k(x)$ must belong to $V$, i.e.
\[  \hat{u}_k(x) = \chi_k \; u_{k-1}(x)+c_0 {\mathscr J} [\delta_{k-1}(x)] \in V. \]
In other words, $\mathscr J$ should be a mapping from $U$ to $\hat{V}$.
In addition, to obey the ``Rule of Completeness'',   $\sum\limits_{k=0}^{+\infty}u_k(x)$  must contain all base functions $\varphi_m \in S_\infty$, $m=1,2,3,\cdots, +\infty$.  Therefore,  ${\mathscr J}[\delta_k(x)]$  as $k\to \infty$  should contain all elements $\varphi_i$ $(i\geq \mu+1)$ of the set $\hat{S}$, where $\hat{S}$ is defined by (\ref{def:S:hat:general}).

In addition, since the high-order deformation equation (\ref{def:mth}) is linear, the inverse operator $\mathscr J$ must be linear, too.  Besides, to guarantee the uniqueness of the solution, $\mathscr J$ must be injective.   Furthermore, the mapping of the inverse operator $\mathscr J$ should be finite.

Therefore, the inverse operator ${\mathscr J}:U\to \hat{V}$ should be defined according to the following rules:
\begin{enumerate}
  \item[(I)] $\mathscr{J}$ is linear, i.e
  \begin{equation*}
   \forall \alpha,\beta \in \mathbb{R}, \forall x, y\in U, \mathscr{J}(\alpha x+\beta y)=\alpha \mathscr{J}(x) + \beta \mathscr{J}(y);
  \end{equation*}
  \item[(II)] $\mathscr{J}$ is  injective, say, the kernel of $\mathscr{J}$ is $\{\mathbf{0}\}$, i.e
  \begin{equation*}
    \{x | x \in U, \mathscr{J}(x) = \mathbf{0}\}=\{\mathbf{0}\};
  \end{equation*}
  \item[(III)]  $\mathscr{J}[\delta_{m}(x)]$ as $m\to +\infty$ contains all base functions $\varphi_i \in \hat{S}$ $(i\geq \mu+1)$;
  \item[(IV)]  $\mathscr{J}$ is finite, i.e. there exists  such a finite constant $K$   that for any $\varphi\in V$  it holds  \[  \frac{||{\mathscr J}[{\mathscr N}[\varphi]]||}{||\varphi||}  \leq  K.    \]
\end{enumerate}

Therefore, using the new strategy of the HAM, the common solution $u_k(x)$ of the $k$th-order deformation equation (\ref{def:mth}) is the sum of the special solution $\hat{u}_k(x)$ and the primary solution $u^*_k(x)$, expressed by
\begin{equation}\label{u_k}
   u_{k}(x) = \chi_k u_{k-1}(x)+ c_{0} \mathscr{J}[  \delta_{k-1}(x)] + \sum^\mu_{i=1} a_{k,i} \; \varphi_i(x),
\end{equation}
where  the constants $a_{k,i}$ ($1\leq i \leq \mu$) is uniquely determined by the $\mu$ linear boundary conditions (\ref{bc:mth}).    We call this new strategy ``the method of directly defining inverse mapping'' (MDDiM).

This is a new strategy to solve differential equation, since we  completely  {\em neglect} the auxiliary linear operator $\mathscr L$ itself, but {\em directly define} its inverse operator $\mathscr J$ using the above rules I - IV.   In this way, we could overcome the restrictions and limitations of  traditional approaches for differential equations!   So, the MDDiM might open a  new way for solving nonlinear differential equations.

\subsection{Proof of the convergence theorem}

It is generally proved \cite{liaobook1, liaobook2} in the frame of the HAM that, if a series solution given by the HAM is absolutely convergent, it must be one solution of original nonlinear equation under consideration.   Since the above-mentioned ``method of direct defining inverse mapping'' (MDDiM) is based on the HAM, one can prove the convergence theorem in a rather similar way.

Since ${\mathscr J}: U\to \hat{V}$ is injective and linear, its inverse operator ${\mathscr L}: \hat{V} \to U$ certainly exists (although we do not know its explicit form) and linear, say,
\begin{enumerate}
  \item[(i)] $\mathscr{L}$ is linear, i.e
  \begin{equation*}
   \forall \alpha,\beta \in \mathbb{R}, \forall x, y\in V, \mathscr{L}(\alpha x+\beta y)=\alpha \mathscr{L}(x) + \beta \mathscr{L}(y);
  \end{equation*}
  \item[(ii)] the composition map $\mathscr{L}\circ \mathscr{J}$ is the identity in $U$, i.e
\[  \forall x\in U,  \mathscr{L}\circ \mathscr{J}[x] = x;\]
  \item[(iii)] ${\mathscr L}[0] = 0$, since ${\mathscr L}$ is injective from $\hat{V} \to U$,
\end{enumerate}
where $V, \hat{V}, U$ are defined by (\ref{def:set:V}), (\ref{def:set:V:hat}) and (\ref{def:set:U}), respectively.

Besides, recall that
\[   u_k^* = \sum_{i=1}^{\mu} a_{k,i} \; \varphi_i  \in V^*  \]
is defined as the primary solution, where $\varphi_i \in S^*$.  Thus,
\[  {\mathscr L}\left[ \sum_{i=1}^{\mu} a_{k,i} \; \varphi_i \right] = 0,   \]
so that  it holds
\begin{enumerate}
\item[(iv)] $ \forall x \in V^*, \; {\mathscr L}[x] = 0.$
\end{enumerate}
In this way,  the linear operator ${\mathscr L}: V \to U$ is well defined.

Here, a proof of the convergence-theorem in \S~2  is given below.

\begin{proof}
Due to (\ref{def:delta-mapping}), it holds  using (i), (ii) and (iv) that
\begin{eqnarray}
{\mathscr L} [u_k] &=&{\mathscr L} \left\{\chi_k \; u_{k-1} + c_0 {\mathscr J}[\delta_{k-1}] + \sum_{i=1}^{\mu}a_{k,i}\; \varphi_i\right\} \nonumber\\
&=& \chi_k \;  {\mathscr L}[ u_{k-1}] + c_0 \delta_{k-1},\nonumber
\end{eqnarray}
since ${\mathscr L} \circ {\mathscr J}[x]=x, \forall x\in U$ and ${\mathscr L}[x]=0, \forall x\in V^*$.  Taking the sum of the above equation from $k=1$ to $+\infty$, we have
\[  \lim_{k\to +\infty} {\mathscr L}[u_k] = c_0 \sum_{n=0}^{+\infty} \delta_n.  \]
If (\ref{def:series-solution}) is absolutely convergent, it holds
\[   \lim_{k\to +\infty} u_k  = 0 \]
which leads to
\begin{eqnarray}
c_0 \sum_{n=0}^{+\infty} \delta_n =\lim_{k\to +\infty} {\mathscr L}[u_k] = {\mathscr L}  [\lim_{k\to +\infty} u_k] =  {\mathscr L}[0] =0.
\end{eqnarray}
Here, the property (iii) of $\mathscr L$ is used.  Furthermore, since $c_0\neq 0$, we have
\begin{equation}
\sum_{n=0}^{+\infty} \delta_n = 0.    \label{sum:delta}
\end{equation}
The Taylor series of
\[  {\mathscr N}\left[ \sum_{n=0}^{+\infty} u_n \; q^n\right] \]
at $q=0$ reads
\[  {\mathscr N}\left[ \sum_{n=0}^{+\infty} u_n \; q^n\right] = \sum_{n=0}^{+\infty} \delta_n \; q^n, \]
which is now convergent to zero at $q=1$, since
\[  {\mathscr N}\left[ \sum_{n=0}^{+\infty} u_n \right] = \sum_{n=0}^{+\infty} \delta_n = 0. \]
Here, (\ref{sum:delta}) is used.    Thus, the series (\ref{def:series-solution}) satisfies the governing equation \[ {\mathscr N}\left[ \sum_{n=0}^{+\infty} u_n\right]=0.\]

Similarly, since  $u_k$ ($k\geq 1$) satisfies the boundary condition (\ref{bc:mth}), we have
\[    {\mathscr B}_i[u_m(x)] = c_i \; \sum_{k=0}^{m-1}\Delta_{i,k}(x), \;\; \mbox{at $x=\alpha_i, \;\; 1\leq i\leq \mu$},   \]
by taking the sum of (\ref{bc:mth}) from $k=1$ to $m$.  As $m\to +\infty$, it becomes
\[    \lim_{m\to +\infty}{\mathscr B}_i[u_m(x)] = c_i \; \sum_{k=0}^{+\infty}\Delta_{i,k}(x), \;\; \mbox{at $x=\alpha_i, \;\; 1\leq i\leq \mu$}.   \]
Similarly, since the solution series  (\ref{def:series-solution}) is absolutely convergent, we have
\[  \lim_{m\to +\infty}{\mathscr B}_i[u_m(x)] = {\mathscr B}_i[\lim_{m\to +\infty}u_m(x)] =  {\mathscr B}_i[0]=0, \]
which leads to
\[   \sum_{k=0}^{+\infty}\Delta_{i,k}(x)=0, \;\; \mbox{at $x=\alpha_i, \;\; 1\leq i\leq \mu$},  \]
since $c_i\neq 0$.   Therefore,  the Maclaurin series
\[   {\mathscr B}_i\left[\sum_{n=0}^{+\infty} u_n(x)\; q^n\right] -\beta_i = \sum_{k=0}^{+\infty}\Delta_{i,k}(x) q^k , \;\;\; 1\leq i\leq \mu,    \]
tends to zero at $q=1$, say,
\[   {\mathscr B}_i\left[\sum_{n=0}^{+\infty} u_n(x)\right] = \beta_i, \;\;\; 1\leq i\leq \mu .   \]
Thus,  the series (\ref{def:series-solution})  satisfies the original boundary condition (\ref{bc:original}), too.

 Therefore,  the solution series (\ref{def:series-solution}) is a solution of (\ref{geq:original}) and (\ref{bc:original}).
 \end{proof}

\section{Illustrative applications}

Three examples are used here to illustrate the validity of the method of directly defining inverse mapping (MDDiM).

\subsection{A nonlinear eigenvalue problem}

First of all, let us consider a nonlinear eigenvalue problem
\begin{equation}\label{geq:exam1}
    {\mathscr N}_1[u,\lambda]=u''(x)+\lambda u(x)+\epsilon u^3(x)=0,
\end{equation}
subject to the boundary conditions
\begin{equation}\label{bc:exam1}
  u(0)=u(1)=0,
\end{equation}
and the normalization condition
\begin{equation}
  \int^1_0 u^2(x) dx=1,  \label{u_n:norm}
\end{equation}
where $\epsilon$ is a physical parameter, the prime denotes differentiation with respect to $x$, respectively.   Here, both of the eigenfunction $u(x)$ and the eigenvalue $\lambda$ are unknown.  This problem has an infinite number of solutions.  Without loss of generality, let us consider here its simplest solution.

According to the odd nonlinearity of Eq.~(\ref{geq:exam1}) and the boundary condition (\ref{bc:exam1}),  $u(x)$ can be expressed by
\begin{equation}\label{form}
u(x)=\sum^{+\infty}_{n=1} a_{n} \sin[(2n-1)\pi x],
\end{equation}
where $a_n$ is a real constant and $n\geq 1$ is an integer.   Regard $\sin(\pi x)$ as the base function of the primary solution.
Then,  we have the following sets
\begin{eqnarray}
V &=& \left\{\left.  \sum_{n=1}^{+\infty} a_n \; \sin[(2n-1)\pi x] \right| a_n \in \mathbb{R} \right\}, \label{def:set:V:exam1}\\
V^* &=& \left\{\left.  a_1 \; \sin[\pi x] \right| a_1 \in \mathbb{R} \right\},\\
\hat{V} &=& \left\{\left.  \sum_{n=2}^{+\infty} a_n \; \sin[(2n-1)\pi x] \right| a_n \in \mathbb{R} \right\}. \label{def:set:V:hat:exam1}
\end{eqnarray}
Note that $V = V^* \cup \hat{V}$.  Thus, it holds $U=V$ for the considered eigenvalue prolem.

The eigenfunction $u(x)$ and eigenvalue $\lambda$ are expressed by
\begin{eqnarray}
u(x)=\sum_{n=0}^{+\infty} u_n(x) \in V, \;\;\; \lambda=\sum_{n=0}^{+\infty} \lambda_n,
\end{eqnarray}
where $u_0(x)\in V^*$ is an initial guess.  Note that $u_n(x) \in V$ for $n\geq 1$ and $\lambda_n\in \mathbb{R}$ for $n\geq 0$.

Consider the Maclaurin series
\[  {\mathscr N}_1 \left[ \sum_{n=0}^{+\infty} u_n(x) q^n,\sum_{n=0}^{+\infty}\lambda_n q^n\right] = \sum_{m=0}^{+\infty}\delta_m(x) \; q^m    \]
about $q$,  where
\begin{eqnarray}
  \delta_{m}(x) &=& u''_{m}(x)+\sum^{m}_{i=0}\lambda_{i}u_{m-i}(x)+\epsilon\sum^{m}_{i=0}u_{m-i}(x)\sum^i_{j=0}u_{j}(x)u_{i-j}(x).\label{def:delta:exam1}
\end{eqnarray}
Obviously,   $\delta_m(x) \in V=U$, where the set $V$ is defined by (\ref{def:set:V:exam1}).

In the frame of the MDDiM, we have
\begin{equation}
u_m(x) = \chi_m \; u_{m-1}(x) + c_0 {\mathscr J}_\alpha[\delta_{m-1}(x)] + a_{m,1} \; \sin(\pi x),  \label{u[m]:exam1}
\end{equation}
where $c_0$ is ``the convergence-control parameter'' whose value we have great freedom to choose, $a_{m,1}$ is a constant to be determined by the normalization condition (\ref{u_n:norm}), ${\mathscr J}_\alpha: V  \to \hat{V}$ is an inverse mapping {\em directly defined} here by
\begin{equation}
\mathscr{J}_\alpha \left\{ \sin[(2m-1) \pi x] \right\}= -\frac{\sin[(2m-1)\pi x]}{2(m-1)(2m+1+\alpha)\pi^2}, \label{def:J:exam1:a}
\end{equation}
where $m > 1$ is an integer and $\alpha>0$ is an auxiliary parameter to be chosen.  Note that different values of $\alpha$ correspond to different inverse mappings.   So,  we actually {\em define} a family of inverse mappings ${\mathscr J}_\alpha$.  According to the Rule (IV),  the inverse mapping ${\mathscr J}_\alpha$ must be finite.  However, ${\mathscr J}_\alpha[\sin(\pi x)]$ tends to infinity.  To avoid this, the term $\sin(\pi x)$ must disappear from $\delta_m(x)$ for $m\geq 0$, say, its coefficient must be zero.  This just provides us an algebraic equation to determine the unknown $\lambda_m$.

Note that the boundary condition $u(0)=u(1)=0$ is automatically satisfied, since $u(x)\in V$, where $V$ is defined by (\ref{def:set:V:exam1}).   Considering the normalization condition (\ref{u_n:norm}), we choose the initial guess $u_0(x) = \sqrt{2}\sin(\pi x)$, since $\sin(\pi x)\in V^*$ is the base function for the primary solution.
 Then, it is straightforward to gain $\delta_{0}(x)$ defined by (\ref{def:delta:exam1}).  Enforcing the coefficient of $\sin(\pi x)$ in $\delta_{0}(x)$ to be zero gives an algebraic equation of $\lambda_{0}$, from which we gain $\lambda_{0}$.  Then, using (\ref{u[m]:exam1}) and the definition (\ref{def:J:exam1:a}) of ${\mathscr J}_\alpha$, we gain $u_1(x)$, whose unknown coefficient $a_{1,1}$ is determined by the normalization condition (\ref{u_n:norm}), i.e.
\begin{equation}
  \int^1_0  \left[ \sum^1_{n=0} u_{n}(x)\right]^2 dx =1. \label{def:norm:1st}
\end{equation}
 In this way, we can gain $\lambda_{0}$, $u_{1}(x)$, $\lambda_{1}$, $u_{2}(x)$, and so on, successively, {\em without} calculating any inverse operators!

In summary, if  $u_0, u_1, \cdots, u_{m-1}$ and $\lambda_0,\lambda_1, \cdots,\lambda_{m-2}$ are known,  it is straightforward to gain $\delta_{m-1}(x)$ defined by (\ref{def:delta:exam1}).  Enforcing the coefficient of $\sin(\pi x)$ in $\delta_{m-1}(x)$ to be zero gives an algebraic equation of $\lambda_{m-1}$, from which we know $\lambda_{m-1}$.  Then, using (\ref{u[m]:exam1}) and the definition (\ref{def:J:exam1:a}) of ${\mathscr J}_\alpha$, we gain $u_m(x)$, whose unknown coefficient $a_{m,1}$ is determined by the normalization condition (\ref{u_n:norm}), i.e.
\begin{equation}
  \int^1_0  \left[ \sum^m_{n=0} u_{n}(x)\right]^2 dx =1. \label{def:norm:mth}
\end{equation}
 In this way, we can gain the series of the eigenvalue $\lambda$ and the eigenfunction $u(x)$, {\em without} calculating any inverse operators.

 To measure the accuracy of the $m$th-order approximation
\[   \bar{u}(x) = \sum_{n=0}^{m}u_n(x), \;\;  \bar{\lambda} = \sum_{n=0}^{m-1}\lambda_n,   \]
we consider the squared residual error
\begin{equation}
  {\mathscr E}_m=\int^1_0 \left\{ \mathcal{N}_1[\bar{u}(x),\bar{\lambda}] \right\}^2 \; dx. \label{def:E[m]}
\end{equation}
Note that the two boundary conditions (\ref{bc:exam1}) are automatically satisfied, and the normalization condition (\ref{u_n:norm}) is also satisfied.
Therefore,  the smaller the squared residual error ${\mathscr E}_m$, the more accurate the $m$th-order approximation $\bar{u}$ and  $\bar{\lambda}$.

It should be emphasized that, unlike perturbation techniques and other traditional methods,  the MDDiM contains an auxiliary parameter $c_0$, called the convergence-control parameter, which  provides us a convenient way to guarantee the convergence of solution series.  For given $\alpha$, the $m$th-order approximation $\bar{u}$ and $\bar{\lambda}$ contain $c_0$.  So does the corresponding residual error square ${\mathscr E}_m$.  Obviously, the optimal value of $c_0$ is determined by the minimum of ${\mathscr E}_m$.   In this example, we use the optimal value $c_0$ gained at the 3rd order of approximation.

Without loss of generality, let us consider the case of $\alpha=2$.  Using the optimal convergence-control parameter $c_0=-5/8$ obtained by the minimum of ${\mathscr E}_3$, we gain a convergent series solution, with ${\mathscr E}_m$ decreasing to $7.5\times10^{-39}$ at the 50th-order of approximation (i.e. $m=50$), as shown in Table~\ref{table:exam1:error}.   This illustrates the validity of the MDDiM.

\begin{table}[ht]
\tabcolsep 0pt \caption{The residual error square ${\mathscr E}_m$ of Eq. (\ref{geq:exam1}) and the relative error of the corresponding eigenvalue $\lambda/\pi^2$ by means of  $\alpha=2$ with the optimal convergence-control parameter  $c_0=-5/8$.}\label{table:exam1:error}
\vspace*{-12pt}
\begin{center}
\def\temptablewidth{1\textwidth}
{\rule{\temptablewidth}{1pt}}
\begin{tabular*}{\temptablewidth}{@{\extracolsep{\fill}}cccc}
$m$, order of approx.  &  ${\mathscr E}_m$ & relative error of $\lambda/\pi^2  \; \; (\%)$ \\
\hline
10&      $7.5\times10^{-6}$   & $7.0\times10^{-4}$\\
20&      $7.5\times10^{-14}$  & $4.8\times10^{-8}$\\
30&      $7.5\times10^{-22}$  & $4.8\times10^{-12}$\\
40&      $1.3\times10^{-30}$  & $5.6\times10^{-16}$\\
50&      $7.5\times10^{-39}$  & $7.0\times10^{-20}$\\
\end{tabular*}
{\rule{\temptablewidth}{1pt}}
\end{center}
\end{table}

\begin{figure}[ht]
\centering
  \includegraphics[width=7cm]{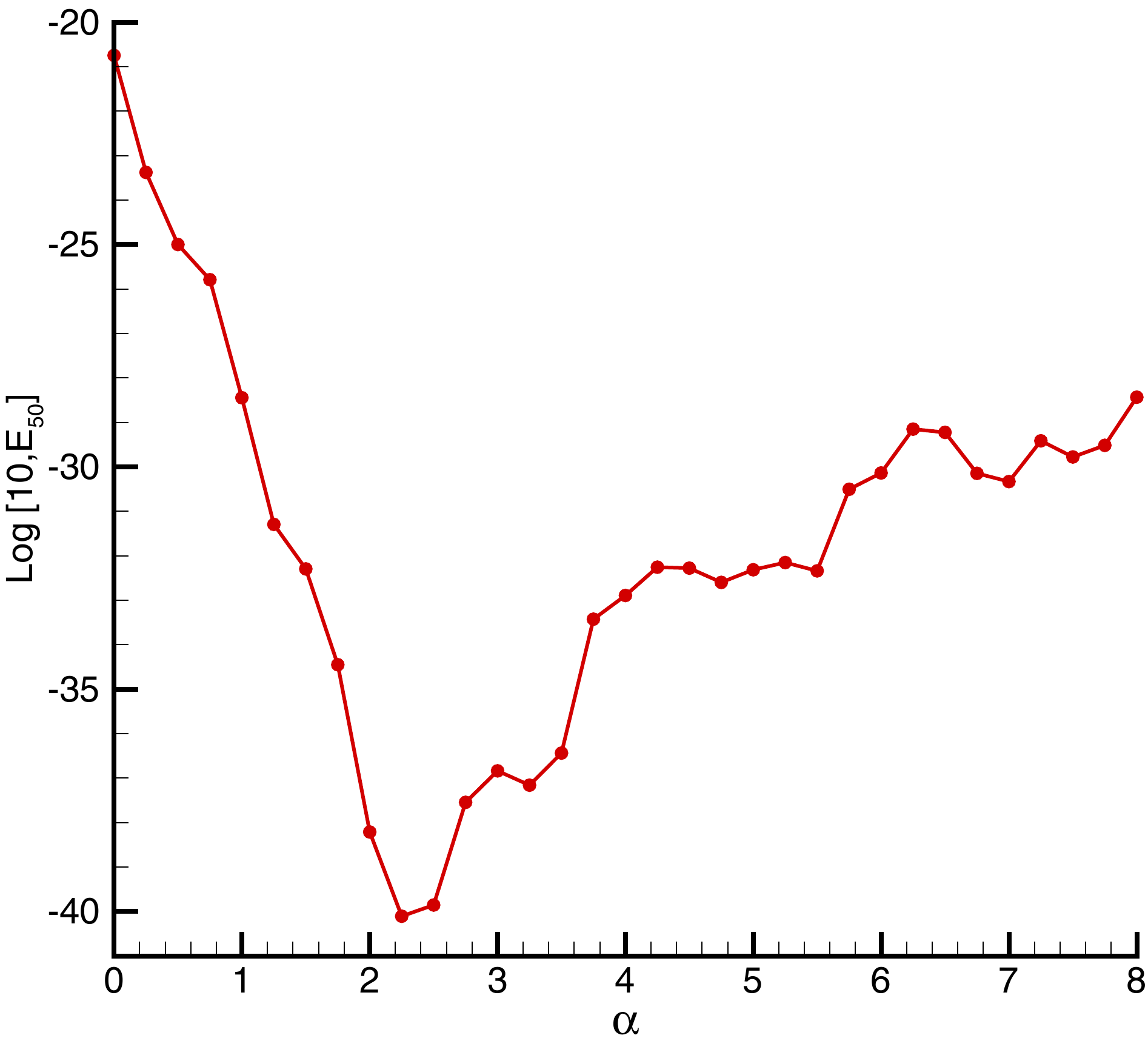}\\
  \renewcommand{\figurename}{Fig.}
  \caption{The residual error square ${\mathscr E}_{50}$ of Eq. (\ref{geq:exam1})  versus the different values $\alpha$ of $\mathscr{J}_{\alpha}$  defined by (\ref{def:J:exam1:a}). }
\label{fig:exam1:error-alpha}
\end{figure}

\begin{figure}[!ht]
\centering
  \includegraphics[width=7cm]{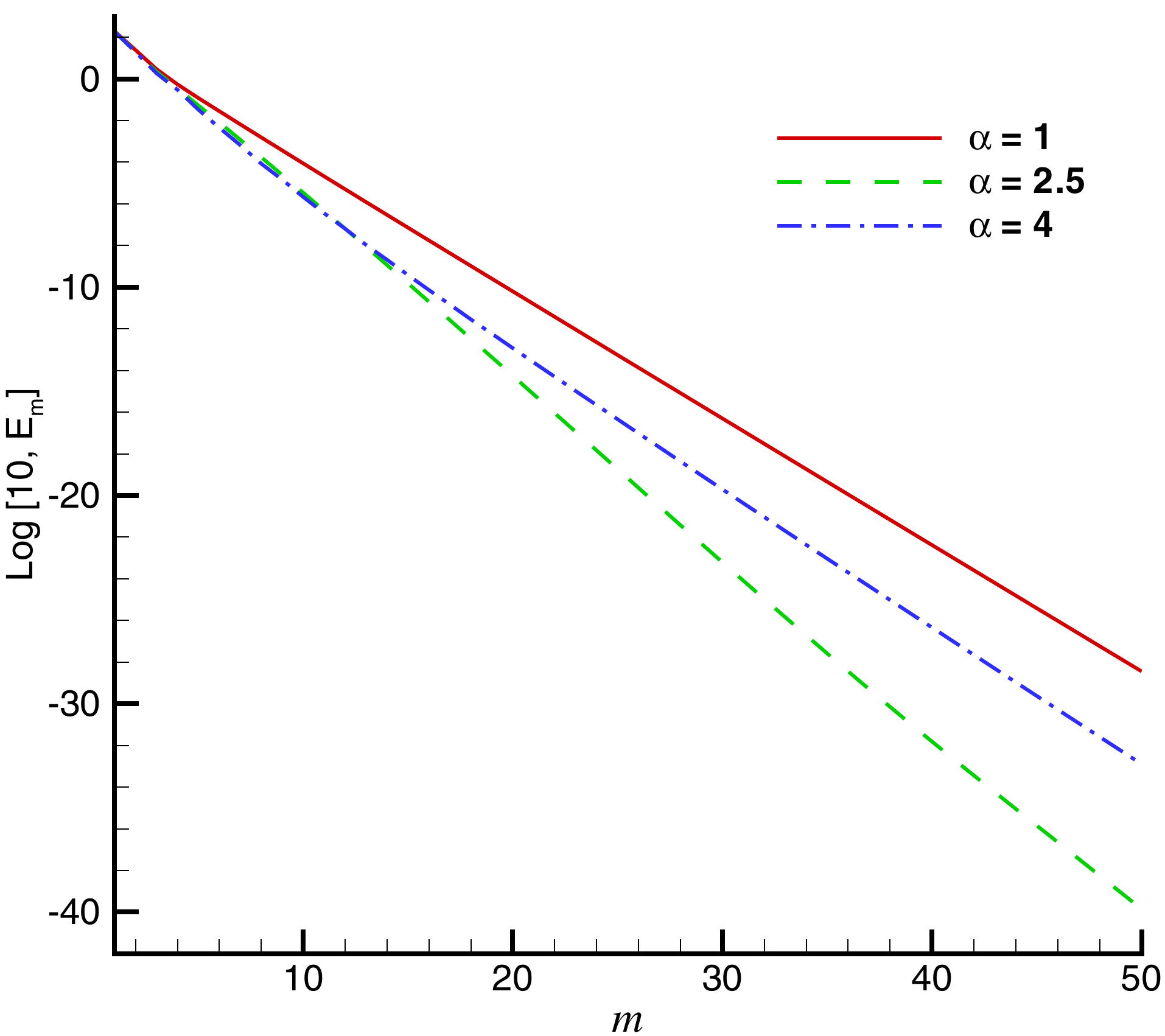}\\
  \renewcommand{\figurename}{Fig.}
  \caption{The residual error square ${\mathscr E}_{m}$ of (\ref{geq:exam1}) versus $m$ (the order of approximation) for the different $\delta$-mapping $\mathscr{J}_{\alpha}$. Solid line: $\alpha=1$ with the optimal value $c_0=-1/2$; Dash-dotted line: $\alpha=4$ with the optimal value $c_0=-11/13$; Dashed line: $\alpha=2.5$ with the optimal value $c_0=-2/3$.}
\label{fig:exam1:error-m}
\end{figure}

Note that we {\em directly define} the inverse mapping (\ref{def:J:exam1:a}) by introducing an auxiliary parameter $\alpha$.
It is found that  we can gain the convergent series solution for {\em any} values of $\alpha\in(0,8)$, as shown in Figs.~\ref{fig:exam1:error-alpha} and \ref{fig:exam1:error-m}, and besides  $\alpha \approx 2.2$ gives the  fastest convergent series.    This  further  illustrates that we {\em indeed} have large freedom and great flexibility to {\em directly define} the inverse mapping ${\mathscr J}_\alpha$.    To confirm this viewpoint, we further consider a more general inverse mapping
\begin{equation}\label{def:J:exam1:b}
  \mathscr{J}_{\beta,\gamma}[\sin (m \pi x)]=\frac{\sin (m \pi x)}{(1-m)(\sqrt{m}+\beta)(\sqrt{m}+\gamma)\pi^2},
\end{equation}
where $m=2k-1$ with $k > 1$.  Using the above inverse mapping with {\em any} values of $\beta \in (0,4)$ and $\gamma\in (0,4)$, we also successfully obtain convergent series solution by means of the corresponding optimal convergence-control parameter $c_0$.  All of these indicate that we {\em indeed} have rather large freedom and great flexibility to {\em directly define} the inverse mapping ${\mathscr J}$ so as to gain the  convergent eigenfunction $u(x)$ and eigenvalue $\lambda$ of Eqs.~(\ref{geq:exam1}) and (\ref{bc:exam1}).

When $\alpha=1$, the corresponding auxiliary operator ${\mathscr L}$ of the inverse mapping  ${\mathscr J}_\alpha$ can be explicitly defined in a differential form, and the considered problem was solved by means of the normal HAM, as mentioned in \S~8 of Liao's book \cite{liaobook1}.  However, as shown in Fig.~\ref{fig:exam1:error-m}, the series given by the MDDiM (when $\alpha=4$ or $\alpha=2.5$)  converge faster even than that given by the normal HAM (corresponding to $\alpha=1$).   It should be emphasized that, in most cases, the two families  (\ref{def:J:exam1:a})  and (\ref{def:J:exam1:b})  of the inverse mapping ${\mathscr J}$ (and its corresponding auxiliary linear operator ${\mathscr L}$) can {\em not} be explicitly defined in a differential form.     The key point is that it is {\em unnecessary} to calculate the auxiliary linear operator $\mathscr L$ at all.   This is more important, since it saves a lots of CPU times and money.    Therefore,  we indeed can {\em directly define} the inverse mapping ${\mathscr J}$ in a more general way.  In other words, the MDDiM is more general than traditional methods that are based on differential operators.  This is the reason why the MDDiM  can give faster convergent series solution in many cases, as shown in this example.   Thus, the MDDiM is {\em fundamentally} different from the traditional methods for differential equations that often spend lots of CPU time to calculate inverse operators.

\subsection{Blasius flow}

Secondly, let us consider the Blasius boundary-layer flow, governed by
\begin{equation}
f'''(\eta) +\frac{1}{2}f(\eta) f'(\eta) = 0, \hspace{0.5cm} f(0)=f'(0)=0, \;\; f'(+\infty) = 1.\label{geq:Blasius}
\end{equation}
Write $ f(\eta) = F(z) + \eta, z=\lambda\; \eta$, where $\lambda>0$ is a constant to be chosen  later.  Then, Eq. (\ref{geq:Blasius}) becomes
\begin{equation}
 {\mathscr N}_2[F] = F''' + \frac{1}{2\lambda^2}(z+\lambda F)F'' = 0,  \label{geq:original:exam3}
\end{equation}
subject to the boundary conditions
\begin{equation}
F(0)=0, F'(0)=-\frac{1}{\lambda},  F'(+\infty)=0,  \label{bc:original:exam3}
\end{equation}
where the prime denotes the derivative with respect to $z$.

In the frame of the MDDiM, we have the solution series
\begin{equation}
F=F_0(z) + \sum_{m=1}^{+\infty} F_m(z),   \label{series:F}
\end{equation}
where $F_0(z)$ is an initial guess satisfying all boundary conditions, and $F_k(z)$ is given by
\begin{equation}
F_m(z) = \chi_{m-1}F_{m-1}(z) + c_0 {\mathscr J}\left[ \delta_{m-1}(z)\right] + F_m^*(z),
\end{equation}
subject to the boundary conditions
\begin{equation}
F_m(0)=F'_m(0)=F'(+\infty)=0,   \label{bc:mthLexam3}
\end{equation}
where $c_0$ is the convergence-control parameter,  $\mathscr J$ is a directly defined inverse mapping,  $F_m^*(z)$ is the primary solution,  and
\begin{eqnarray}
\delta_k(z) &=& {\mathscr D}_k  \left\{ {\mathscr N}_2 \left[\sum_{n=0}^{+\infty} F_n(z) \; q^n  \right] \right\} \nonumber\\
&=& F_k'''(z) + \frac{z}{2\lambda^2}F_k''(z)+\frac{1}{2\lambda}\sum_{n=0}^{k}F_{k-n}(z) F''_n(z), \label{def:delta:exam3}
\end{eqnarray}
respectively.

\begin{figure}[!htb]
\centering
  \includegraphics[width=9cm]{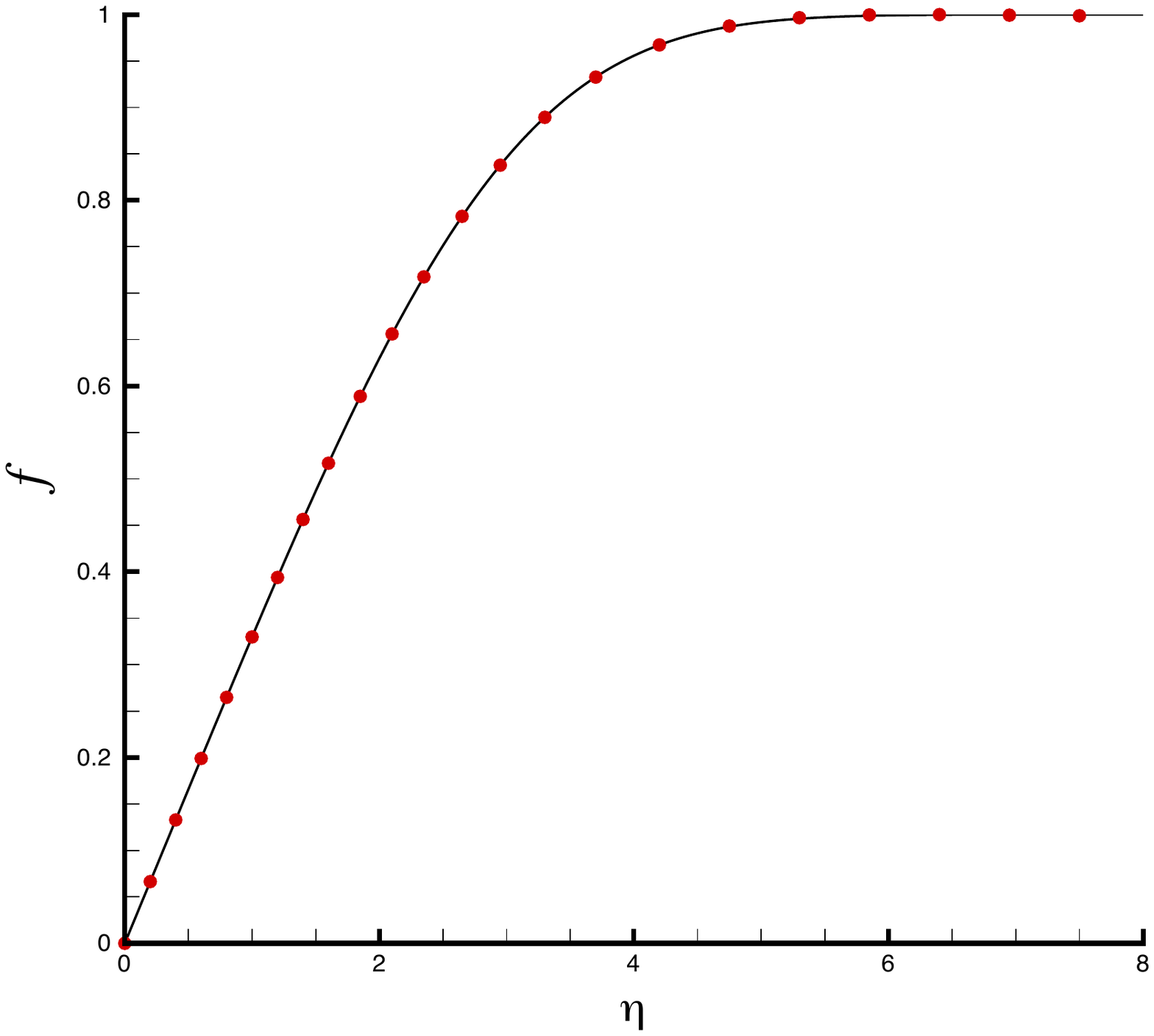}\\
  \renewcommand{\figurename}{Fig.}
  \caption{Comparison of $f'(\eta)$ between the numerical result and the 30th-order approximation given by means of the MDDiM using the directly defined inverse mapping (\ref{def:J:exam3}) with $A_0=1/(3\pi), A_1=\pi/30, A_2=\pi/3$ and $\lambda=1/3,  c_0=-9/5$.  Solid line:  numerical result; Symbols: analytic result given by the MDDiM. }
\label{fig:exam3:Blasius}
\end{figure}

According to (\ref{bc:mthLexam3}),  $F'(z)$ tends to zero at infinity.   So, we define the sets
\begin{eqnarray}
V &=& \left\{\left.  \sum_{n=0}^{+\infty} \frac{a_n}{(1+z)^n}  \right| a_n \in \mathbb{R} \right\}, \label{def:set:V:exam3}\\
\hat{V} &=& \left\{\left.  \sum_{n=2}^{+\infty} \frac{a_n}{(1+z)^n }  \right| a_n \in \mathbb{R} \right\} = U, \label{def:set:V:hat:exam3}\\
V^* &=& \left\{\left.  \sum_{n=0}^{1} \frac{a_n}{(1+z)^n }  \right| a_n \in \mathbb{R} \right\}. \label{def:set:V:star:exam3}\\
\end{eqnarray}
Note that $V=\hat{V} \cup V^*$.
Obviously, $F(z) \in V$ and $\delta_k(z) \in \hat{V} = U$.  Obviously, it is straight forward to choose the initial guess
\begin{equation}
F_0(z)=\frac{1}{\lambda}\left( \frac{1}{1+z}-1\right) \in V^*,
\end{equation}
which satisfies all boundary conditions (\ref{bc:original:exam3}).   Besides,  according to (\ref{def:delta:exam3}), $\delta_k(z)$ can be expressed by
\[   \delta_k(z) =  \sum_{m=2}^{+\infty} \frac{a_{k,m}}{(1+z)^m} \in U,   \]
where $a_{k,m}$ is a real coefficient.

\begin{table}[tb]
\tabcolsep 0pt
\caption{Approximations of Blasius boundary-layer flows by means of the MDDiM  using the directly defined inverse mapping  (\ref{def:J:exam3}) with $A_0=1/(3\pi)$, $A_1=\pi/30, A_2=\pi/3$ and $\lambda=1/3$, $c_0=-9/5$.}\label{table:exam3:error}
\vspace*{-12pt}
\begin{center}
\def\temptablewidth{1\textwidth}
{\rule{\temptablewidth}{1pt}}
\begin{tabular*}{\temptablewidth}{@{\extracolsep{\fill}}cccc}
$m$, order of approx.  &   $f''(0)$  &  ${\mathscr E}_m$  \\
\hline
10&      0.34354   & $2.3\times10^{-3}$\\
20&      0.33362   & $8.5\times10^{-5}$\\
30&      0.33206   & $2.3\times10^{-6}$\\
40&      0.33213   & $9.8\times10^{-8}$\\
50&      0.33207   & $2.0\times10^{-8}$\\
60&      0.33203   & $3.2\times10^{-9}$ \\
70&      0.33207   & $7.0\times10^{-10}$ \\
80&      0.33207   & $3.2\times10^{-10}$ \\
90&      0.33205   & $1.1\times10^{-10}$ \\
100&     0.33205   & $2.2\times10^{-11}$ \\
\end{tabular*}
{\rule{\temptablewidth}{1pt}}
\end{center}
\end{table}

  In the frame of the MDDiM, we directly define such an inverse mapping ${\mathscr J}: U \to \hat{V}$ that
\begin{eqnarray}
{\mathscr J}\left[  (1+z)^m \right] = \frac{(1+z)^m}{m^3+A_2 m^2 + A_1 m+A_0}, \;\;\; m\leq -2,  \label{def:J:exam3}
\end{eqnarray}
where $A_0, A_1$ and $A_2$ are constants to be chosen.   Its special solution reads
\begin{equation}
\hat{F}_m =  \chi_m \; F_{m-1} + c_0 {\mathscr J}\left[\delta_{m-1}\right].
\end{equation}
and  the primary solution is
\begin{equation}
F_m^* = a_{m,0} + \frac{a_{m,1}}{1+z}  \in V^*,
\end{equation}
where $a_{m,0}$ and $a_{m,1}$ are real coefficients.  Thus, we have the solution
\begin{equation}
F_m(z) = \hat{F}_m + F_m^* = \chi_m \; F_{m-1} + c_0 {\mathscr J}\left[{\delta}_{m-1}(z)\right] + a_{m,0} + \frac{a_{m,1}}{(1+z)},
\end{equation}
where $a_{m,0}$ and $a_{m,1}$ are determined by   $F_m(0) = F'_m(0) = 0$ of the boundary conditions (\ref{bc:mthLexam3}), since $F_m'(+\infty)$ is automatically satisfied.

In the frame of the MDDiM,  the ``convergence-control parameter'' $c_0$ provides us a convenient way to guarantee the convergence of solution series.  For properly chosen parameters $A_2, A_1, A_0$ of the inverse mapping $\mathscr J$ defined by (\ref{def:J:exam3}), one can choose an optimal value of the convergence-control parameter $c_0$ for a fastest convergence of the series (\ref{series:F}).  For example, we can gain the convergent series solution by means of
\[  \lambda=\frac{1}{3}, \;\; A_0=\frac{1}{3\pi}, \;\;  A_1 =\frac{\pi}{30}, \;\; A_2 = \frac{\pi}{3}, \;\; c_0 = -\frac{9}{5},    \]
as shown in Table~\ref{table:exam3:error}.    The corresponding 30th-order approximations agrees well with the numerical ones in the whole interval $\eta\in[0,+\infty)$, as shown in Fig.~\ref{fig:exam3:Blasius}.   It is found that such kind of inverse mapping $\mathscr J$ is {\em not} unique:  one can gain convergent series solution by means of many inverse mappings, such as
\[  \lambda=\frac{1}{3}, \;\; A_0=0, \;\;  A_1 =0, \;\; A_2 = \frac{\pi}{3}, \;\; c_0 = -\frac{3}{2},    \]
or
\[  \lambda=\frac{1}{3}, \;\; A_0=\frac{1}{10}, \;\;  A_1 =\frac{\pi}{12}, \;\; A_2 = \frac{\pi}{3}, \;\; c_0 = -\frac{3}{2},    \]
and so on: all of them give the {\em same} results that converge to the numerical ones!

This example illustrates  that, in the frame of the MDDiM,  there indeed exist many {\em directly defined} inverse mappings $\mathscr J$, which lead to the {\em same} convergent series solutions of Blasius boundary-layer flow, as long as they are properly defined.  The 2nd example shows once again the validity and potential of the MDDiM.

\subsection{Gelfand equation}

Finally, let us consider the two-dimensional Gelfand equation
\begin{equation}
\nabla^2 u + \lambda \; e^u = 0,  \hspace{0.5cm}  x\in[-1,1], y\in[-1,1],   \label{geq:Gelfand:original:u}
\end{equation}
subject to the boundary conditions
\begin{equation}
u(x,\pm1)=f(x,\pm 1), \;\;  u(\pm 1, y) = f(\pm 1,y), \label{bc:Gelfand:original:u}
\end{equation}
where $u(x,y)$ is the unknown eigenfunction, $\lambda$ is the unknown eigenvalue, and $f(x,y)$ is a given smooth even function, respectively.

Define $u(0,0)=A$ and write $u=A+w$, where $A$ is a unknown constant.   The above equations becomes
 \begin{equation}
{\mathscr N}_3[w,\lambda]=\nabla^2 w + \left( \lambda e^A\right)  \; e^w = 0,  \hspace{0.5cm}  x\in[-1,1], y\in[-1,1],   \label{geq:Gelfand:original}
\end{equation}
subject to the boundary conditions
\begin{equation}
w(x,\pm1)=-A+f(x,\pm 1), \;\;  w(\pm 1, y) = -A+f(\pm 1, y),    \label{bc:Gelfand:original}
\end{equation}
with the restriction
\begin{equation}
w(0,0)=0.    \label{bc:Gelfand:original:2}
\end{equation}

Obviously, for a given $A$,  if $w(x,y)$ and $\lambda$ satisfy the governing equation (\ref{geq:Gelfand:original}) and the boundary conditions  (\ref{bc:Gelfand:original}), then all of $w(-x,y)$, $w(x,-y)$ and $w(-x,-y)$ are its solutions, since $f(x,y)$ is an even function.  So, $w(x,y)$ is an even function of $x$ and $y$, and thus can be expressed by
\begin{equation}
w(x,y) = \sum_{m=0}^{+\infty}\sum_{n=0}^{+\infty} a_{m,n} \; x^{2m} \; y^{2n}.
\end{equation}
Define the sets
\begin{equation}
V = U = \left\{ \left. \sum_{m=0}^{+\infty}\sum_{n=0}^{+\infty} a_{m,n} \; x^{2m} \; y^{2n}  \right| a_{m,n} \in \mathbb{R}\right\} \label{def:set:V:Gelfand}
\end{equation}
and
\begin{eqnarray}
\hat{V} &=& \left\{ \left. \sum_{m=1}^{+\infty}\sum_{n=1}^{+\infty} a_{m,n} \; x^{2m} \; y^{2n}  \right| a_{m,n} \in \mathbb{R}\right\}. \label{def:set:V:hat:Gelfand}
\end{eqnarray}

In the frame of the MDDiM, we have the $m$th-order approximation
\[   w \approx w_0(x,y) + \sum_{n=1}^{m} w_n(x,y),  \;\; \lambda \approx \sum_{n=0}^{m} \lambda_n, \]
where $w_0(x,y)$ is the initial guess, and
\begin{equation}
w_n(x,y)  =  \hat{w}_n + w^*_n, \;\;\; n\geq 1,
\end{equation}
in which
\begin{equation}
\hat{w}(x,y) =  \chi_n \; w_{n-1}(x,y) + c_0 \; {\mathscr J}\left[ \delta_{n-1}(x,y) \right]
\end{equation}
is a special solution, $w^*(x,y)$ is a primary  solution, $c_0$ is the ``convergence-control parameter'',  and
\begin{eqnarray}
\delta_k &=&  {\mathscr D}_k \left\{ {\mathscr N}_3\left[\sum_{i=0}^{+\infty}\lambda_i q^i, \sum_{i=0}^{+\infty}u_i\; q^i  \right] \right\}
  \nonumber \\
&=&\nabla^2 w_k + e^A\sum_{i=0}^{k}\lambda_{k-i} G_{k}(x,y),
\end{eqnarray}
with the definition
\[   G_0=e^{u_0}, \;\;   G_k = \sum_{i=0}^{k-1}\left( 1-\frac{i}{k}\right)w_{k-i} G_i, \]
 respectively.  Note that $\delta_k \in U$.  Thus,  in the frame of the  MDDiM, we  {\em directly define}  an inverse mapping ${\mathscr J}$: $U\to\hat{V}$, say,
\begin{equation}
{\mathscr J}\left[x^m \; y^n  \right] = \frac{x^{m+2} \; y^{n+2}}{(m^2+B_1 m + B_0)(n^2+B_1 n + B_0)}, \;\; m\geq 0,n \geq 0,  \label{def:J:Gelfand2D-A}
\end{equation}
where $B_0 > 0$ and $B_1 > 0$ are constants, and $U$ and $\hat{V}$ are defined by (\ref{def:set:V:Gelfand}) and (\ref{def:set:V:hat:Gelfand}), respectively.

For the sake of the completeness,  we have the primary solution $w_n^* \in V^*$, where
 \begin{equation}
V^* = \left\{ \left. \sum_{m=0}^{+\infty} \left(a_{m,0} \; x^{2m} + a_{0,m} \; y^{2n}\right)  \right| a_{m,0},a_{0,m}  \in \mathbb{R}\right\}, \label{def:set:V:star:Gelfand}
\end{equation}
since $w_n \in V$ and $\hat{w}_n \in \hat{V}$.  The primary solution $w_n^*$ is determined by the boundary conditions
\begin{equation}
w_n-\chi_n\; w_{n-1} = c_1 \; \left[  w_{n-1} + (1-\chi_n)[A-f(x,y)]  \right], \mbox{at $x=\pm 1, y=\pm 1$},
\end{equation}
where $c_1$ is the 2nd  ``convergence-control parameter''.   For simplicity, write
\begin{equation}
w_n = \Gamma_n(x,y), \hspace{0.5cm} \mbox{at $x=\pm 1$ and $y=\pm 1$}, \label{bc:Gelfand:mth}
\end{equation}
where
\[ \Gamma_n(x,y) = \chi_n\; w_{n-1}+c_1 \; \left[  w_{n-1} + (1-\chi_n)[A-f(x,y)]  \right]. \]
Substituting $w_n=\hat{w}_n+w^*_n$ into the boundary conditions (\ref{bc:Gelfand:mth}), we have the primary solution
\begin{eqnarray}
w_n^*(x,y) &=& -\hat{w}_n(x,\pm 1) - \hat{w}_n(\pm 1,y) + \hat{w}_n(\pm,\pm 1)  \nonumber\\
&+& \Gamma_n(x,\pm 1) + \Gamma_n(\pm 1,y) -\Gamma_n(\pm 1, \pm 1).
\end{eqnarray}
Finally, we have the solution
\begin{eqnarray}
w_n(x,y) &=&   \hat{w}_n(x,y)  -\hat{w}_n(x,\pm 1) - \hat{w}_n(\pm 1,y) + \hat{w}_n(\pm,\pm 1)  \nonumber\\
&+& \Gamma_n(x,\pm 1) + \Gamma_n(\pm 1,y) -\Gamma_n(\pm 1, \pm 1),
\end{eqnarray}
which satisfies all of the boundary conditions (\ref{bc:Gelfand:mth}).  Up to now, $\lambda_{k-1}$ is unknown. Note that, according to the restriction condition (\ref{bc:Gelfand:original:2}), we have $w_n(0,0)=0$.  This just provides us an algebraic equation for the unknown  $\lambda_{n-1}$.  For simplicity, we choose the initial guess $w_0(x,y)=0$.  Then, using the above approach, we can gain $w_1, \lambda_0$, then $w_2, \lambda_1$, and so on, step by step.

\begin{figure}[!ht]
\centering
  \includegraphics[width=7cm]{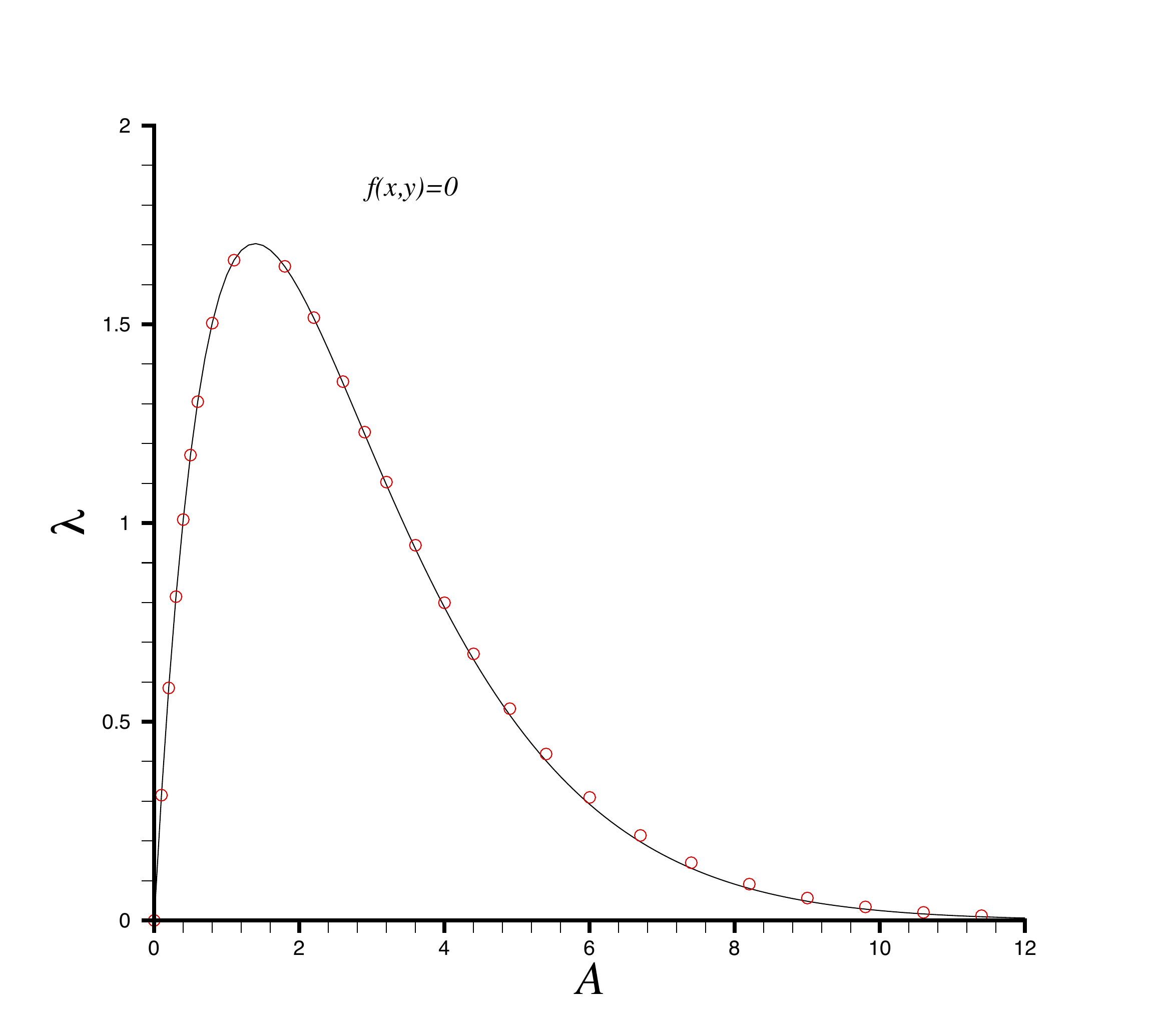}\\
  \renewcommand{\figurename}{Fig.}
  \caption{The eigenvalue of the Gelfand equation in case of $f(x,y)=0$ by means of the MDDiM using  the directly defining inverse mapping  (\ref{def:J:Gelfand2D-A}) with $B_1=\pi, B_0=\pi/2$ and $c_0=3/4, c_1=-3/4$.  Solid line: 20th-order approximation; Symbols: 25th-order approximation.}
\label{fig:exam2:Gelfand2D-1}
\end{figure}

\begin{figure}[!ht]
\centering
  \includegraphics[width=7cm]{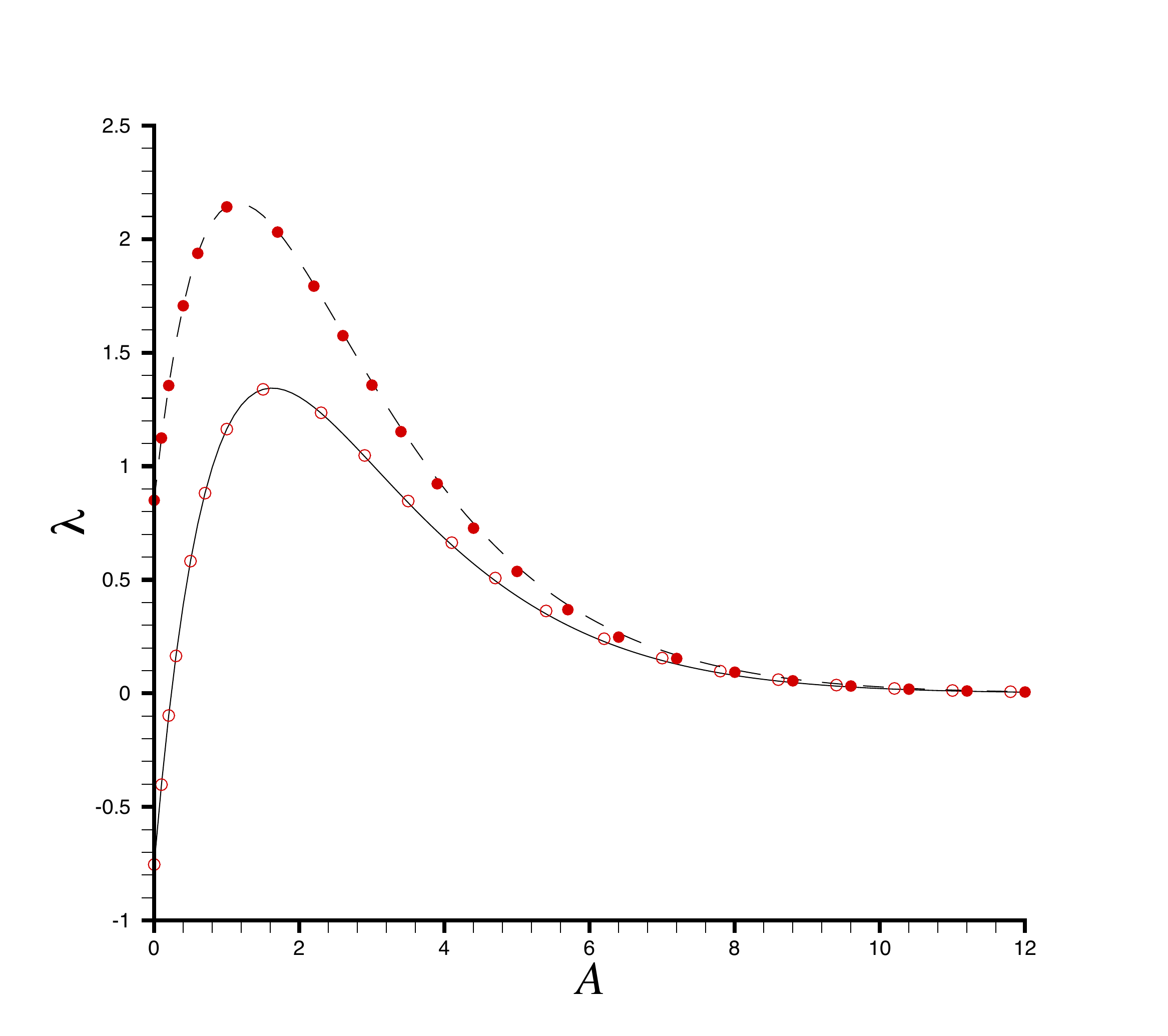}\\
  \renewcommand{\figurename}{Fig.}
  \caption{The eigenvalue of the Gelfand equation in case of $f(x,y)=\pm (1+x^2)(1+y^2)/10$ by means of the MDDiM using  the directly defining inverse mapping  (\ref{def:J:Gelfand2D-A}) with $B_1=\pi, B_0=\pi/2$ and $c_0=1/2, c_1=-1/2$.  Lines: 20th-order approximation; Symbols: 25th-order approximation. Solid line: $f(x,y)=(1+x^2)(1+y^2)/10$; Dashed line: $f(x,y)=-(1+x^2)(1+y^2)/10$. }
\label{fig:exam2:Gelfand2D-2}
\end{figure}

Note that there exist two convergence-control parameters $c_0$ and $c_1$.  Besides, we have great freedom to choose the two  auxiliary parameters $B_0$ and $B_1$ in the directly defined inverse mapping (\ref{def:J:Gelfand2D-A}).    It is found that the convergent series solution can be obtained by means of choosing proper convergence-control parameters $c_0, c_1$ and the two auxiliary parameters $B_0, B_1$ in (\ref{def:J:Gelfand2D-A}).   For example, in case of $f(x,y)=0$,    we gain the good approximation of $u(x,y)$ and $\lambda$ for $A\in[0,12]$ by means of choosing $B_1=\pi, B_0=\pi/2$ and $c_0=3/4, c_1=-3/4$,  as shown in Fig.~\ref{fig:exam2:Gelfand2D-1}.   Such kind of inverse mapping is {\em not} unique: the same convergent result can be obtained by means of choosing $B_1=3, B_0=2$ and $c_0=1, c_1=-1$.  This illustrates that we indeed have large freedom  and great flexibility to directly define the inverse mapping (\ref{def:J:Gelfand2D-A})!

Similarly, in case of
\begin{equation}
f(x,y) = \pm \frac{(1+x^2)(1+y^2)}{10},  \label{def:f:2}
\end{equation}
we gain the good approximation of $u(x,y)$ and $\lambda$ for $A\in[0,12]$ by means of  $B_1=\pi, B_0=\pi/2$ and $c_0=1/2, c_1=-1/2$, as shown in Fig.~\ref{fig:exam2:Gelfand2D-2}.

\begin{figure}[tb]
\centering
  \includegraphics[width=7cm]{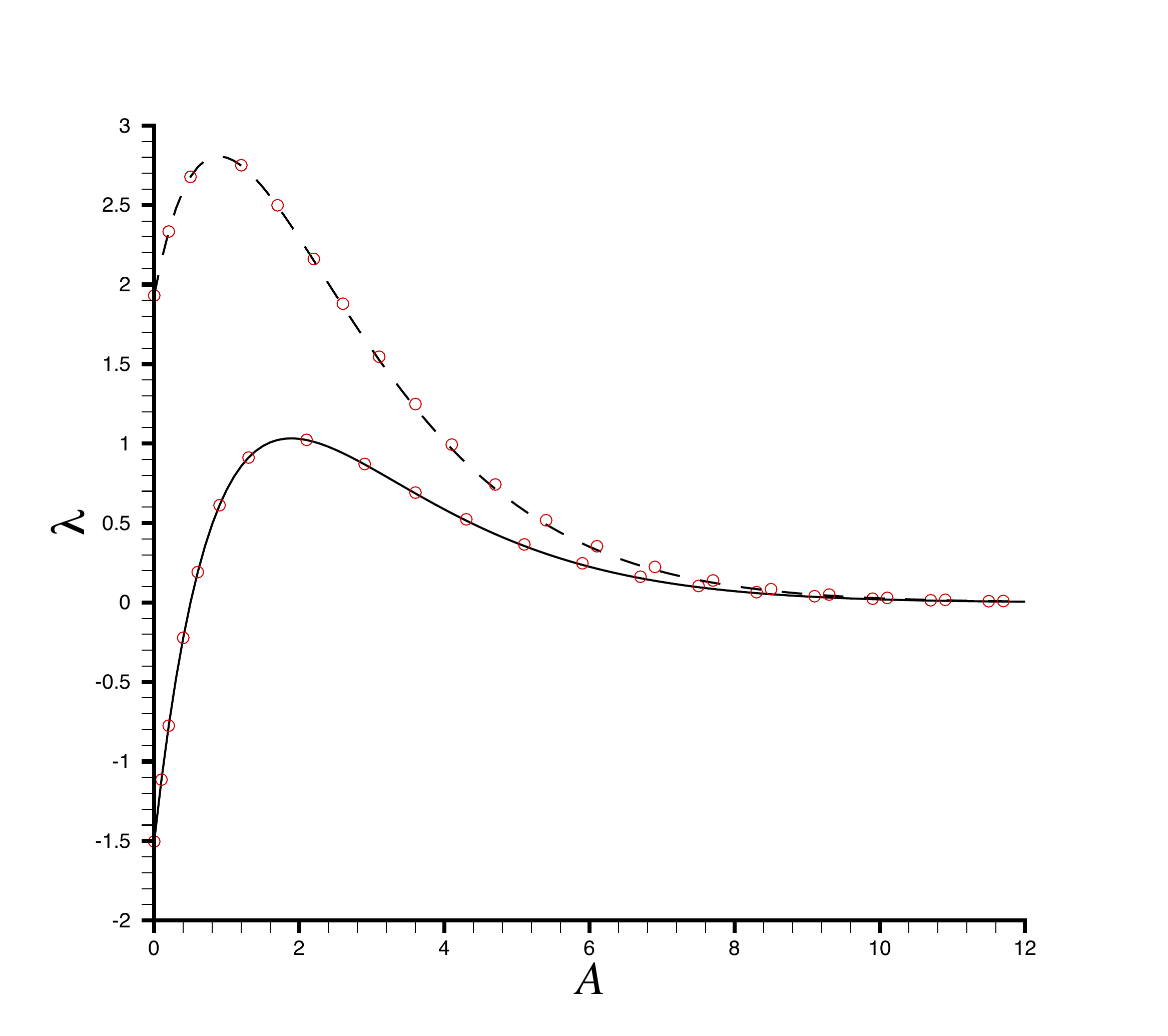}\\
  \renewcommand{\figurename}{Fig.}
  \caption{The eigenvalue of the Gelfand equation in case of $f(x,y)=\pm (x^2-x^2y^2+y^2)/2$ by means of the MDDiM using  the directly defining inverse mapping  (\ref{def:J:Gelfand2D-A}) with $B_1=\pi, B_0=\pi/2$ and $c_0=1/2, c_1=-1/2$.  Lines: 20th-order approximation; Symbols: 25th-order approximation. Solid line: $f(x,y)=(x^2-x^2y^2+y^2)/2$; Dashed line: $f(x,y)=-(x^2-x^2y^2+y^2)/2$. }
\label{fig:exam2:Gelfand2D-3}
\end{figure}

\begin{figure}[!ht]
\centering
  \includegraphics[width=7cm]{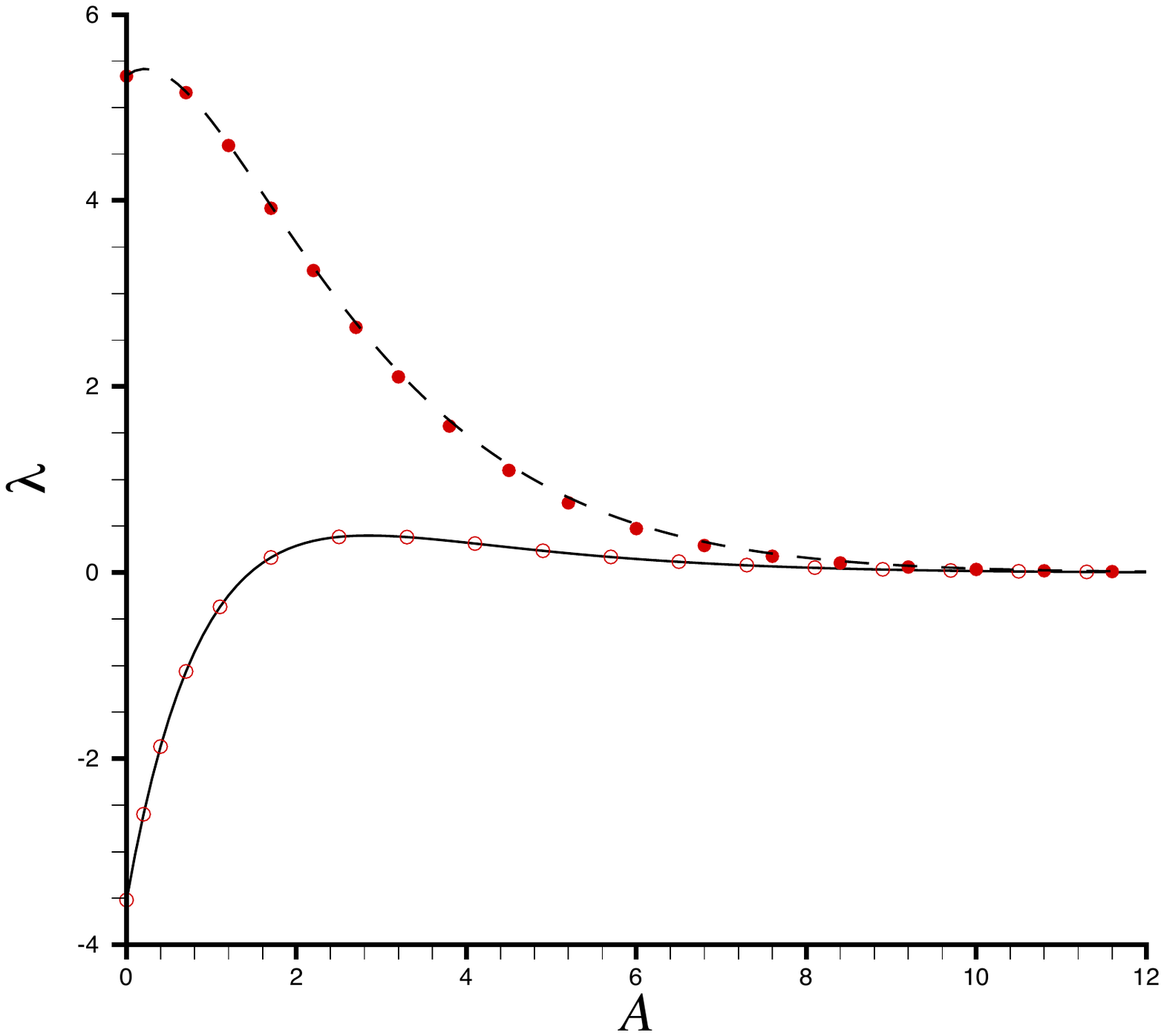}\\
  \renewcommand{\figurename}{Fig.}
  \caption{The eigenvalue of the Gelfand equation in cases of $f(x,y)=\cos x +\cos y$ and $f(x,y)=\cos[\sin(x)]-\exp(y^2)$ by means of the MDDiM using  the directly defining inverse mapping  (\ref{def:J:Gelfand2D-A}) with $B_1=\pi, B_0=\pi/2$ and $c_0=1, c_1=-1$.  Lines: 15th-order approximation; Symbols: 20th-order approximation. Solid line: $f(x,y)=\cos x +\cos y$; Dashed line: $f(x,y)=\cos[\sin(x)]-\exp(y^2)$. }
\label{fig:exam2:Gelfand2D-4}
\end{figure}

In case of
\begin{equation}
f(x,y)=\pm (x^2-x^2y^2+y^2)/2,  \label{def:f:3}
\end{equation}
the good approximation of $u(x,y)$ and $\lambda$ for $A\in[0,12]$ are gained by means of  $B_1=\pi, B_0=\pi/2$ and $c_0=1/2, c_1=-1/2$, as shown in Fig.~\ref{fig:exam2:Gelfand2D-3}.  In case of
\begin{equation}
f(x,y) = \cos (x) + \cos (y), \label{def:f:4}
\end{equation}
the good approximation of $u(x,y)$ and $\lambda$ for $A\in[0,12]$ are gained by means of $B_1=\pi, B_0=\pi/2$ and $c_0=1, c_1=-1$, as shown in Fig.~\ref{fig:exam2:Gelfand2D-4}.  Here, we use
\[    \cos(x) + \cos(y) \approx 2-\frac{1}{2}(x^2+y^2) +\frac{1}{24}(x^4+y^4)-\frac{1}{720}(x^6+y^6),   \]
which is a good approximation for all $x\in[-1,1]$ and $y\in[-1,1]$.   In case of
\begin{equation}
f(x,y) = \cos [\sin(x)] - \exp(y^2), \label{def:f:5}
\end{equation}
the good approximation of $u(x,y)$ and $\lambda$ for $A\in[0,12]$ are gained by means of  $B_1=\pi, B_0=\pi/2$ and $c_0=3/4, c_1=-3/4$, as shown in Fig.~\ref{fig:exam2:Gelfand2D-4}.  Here, we use
\begin{eqnarray}
&& \cos [\sin(x)] - \exp(y^2) \nonumber \\
&\approx&
 -\frac{1}{2}x^2+\frac{5}{24}x^4-\frac{37}{720}x^6+\frac{457}{40320}x^8-\frac{389}{172800}x^{10}   \nonumber\\
 &&- \left(y^2+\frac{1}{2}y^4+\frac{1}{6}y^6+\frac{1}{24}y^8+\frac{1}{120} y^{10}  \right),
\end{eqnarray}
which is a good approximation for all $x\in[-1,1]$ and $y\in[-1,1]$.   Thus, by means of the MDDiM,  the two-dimensional Gelfand equation (\ref{geq:Gelfand:original:u}) and (\ref{bc:Gelfand:original:u}) with rather complicated even function $f(x,y)$  can be easily solved in a straight-forward way.   Note that the inverse mapping $\mathscr J$ (which leads to convergent results) is {\em not} unique in all of these cases.

Finally, it should be mentioned that, when $B_1=3$ and $B_0=2$, the directly defined inverse mapping (\ref{def:J:Gelfand2D-A}) can be expressed in a differential form
\begin{equation}
{\mathscr L} u ={\mathscr J}^{-1} u = \frac{\partial^4 u}{\partial x^2 \partial y^2},
\end{equation}
which is used by Liao and Tan \cite{Liao2007} in the frame of the normal HAM.  As mentioned in \cite{Liao2007}, the original 2nd-order Gelfand equation is transferred into an infinite number of 4th-order linear differential equations governed by an auxiliary linear operator $\mathscr L$ defined above.  This is very difficult  to understand in the frame of the traditional methods for differential equations, which often transfer a $n$th-order differential equation to some sub-equations but only with the {\em same} order.    However, in the frame of the MDDiM, it is easy and straight-forward to understand it, since the MDDiM is based on a mapping that is {\em more general} than a  differential operator.   Especially, it should be emphasized that, when $B_2=\pi$ and $B_0=\pi/2$ (as we used in this paper),  the auxiliary linear operator ${\mathscr L}={\mathscr J}^{-1}$ can {\em not} be expressed in a differential form!   Fortunately, we now need {\em not} consider the auxiliary linear operator $\mathscr L$ at all,  mainly because the MDDiM is based on the {\em directly defined} inverse mapping $\mathscr J$, {\em without} considering its original auxiliary operator $\mathscr L$.  This opens a  new, more general way to solve nonlinear differential equations, which is fundamentally different from the traditional methods.

All of these examples illustrate the validity of the MDDiM, and especially the great freedom and large flexibility of directly defining the inverse mapping $\mathscr J$ for various types of nonlinear problems.

\section{Concluding remarks}

In scientific computation,  it is time-consuming to calculate inverse operators of a differential equation.   Can we solve a nonlinear differential equation without calculating any inverse operators?

The answer is positive:  we can indeed solve nonlinear differential equations by directly defining an inverse mapping $\mathscr J$, as described in this article.  In this work, the ``method of directly defining inverse mapping'' (MDDiM) is proposed based on the homotopy analysis method (HAM)  \cite{liaoPhd, liaobook1, liaobook2, liaobook3}, a widely used analytic approximation technique for highly nonlinear problems.   By means of the MDDiM, one indeed can solve a nonlinear differential equation {\em without} searching for any inverse operators at all, as illustrated in this paper.  From this viewpoint,  the MDDiM is fundamentally different from the traditional ones, which often spend lots of time to calculate inverse operators.

To simplify the use of the MDDiM, some rules are given to guide how to directly define an inverse mapping $\mathscr J$.  Besides, a convergence theorem is proved, which guarantees that a convergent series solution given by the MDDiM must be one solution of problems under consideration.  In addition, three examples are used to illustrate the validity and potential of the MDDiM.

The MDDiM can be regarded as a generalization of the HAM and other traditional methods: it  directly uses mappings between base functions,  instead of differential operators.   Note that  mapping is more general than differential operator.  So, theoretically speaking,  the MDDiM is more general than the normal HAM and other traditional methods which are based on differential operators.

In the frame of the normal HAM,  the 2nd-order two-dimensional Gelfand equation were replaced by an infinite number of the 4th-order (two-dimensional) linear differential equations, as shown by Liao and Tan \cite{Liao2007} who gained accurate approximations with good agreement to numerical ones.    However, this is very difficult to understand in the frame of the traditional methods for differential equations.   But, from the viewpoint of the MDDiM,  it is easy and straight-forward to understand, since the MDDiM gives up the concept of ``differential operator'' at all: it is based on directly defining inverse {\em mapping} that is a concept more general than ``differential operator''.

Note that many differential equations have their equivalent form in integral.   This suggests that many integral equations can be solved by means of MDDiM.  Although the three examples used in this paper are boundary-value problems,  the MDDiM should be also valid for some initial problems whose solutions are not chaotic.

In summary, the MDDiM might bring us a new, more general way to solve nonlinear differential equations, if base functions and inverse mapping are properly chosen.  Without doubt, the MDDiM is at its very beginning, and thus further theoretical researches and  more applications are certainly needed in future.

\section*{Acknowledgment}

This work is supported by National Natural Science Foundation of China (Approval No. 11272209) and State Key Laboratory of Ocean Engineering (Approval No.GKZD010065).


\appendix
\section{The properties of the homotopy-derivative $\mathscr{D}_m$}

For two series
\[   \phi(x;q) = \sum_{k=0}^{+\infty} u_k(x) \; q^k,  \hspace{1.0cm} \psi(x;q) = \sum_{k=0}^{+\infty} w_k(x) \; q^k,    \]
where $\phi(x;q)$ and $\psi(x;q)$ are analytic in $q \in [0,a]$, it holds  for integer $m\geq 0$ that
\begin{eqnarray}
\mathscr{D}_m [\phi] &=& u_m,  \\
\mathscr{D}_m [q^k  \phi] &=& \mathscr{D}_{m-k}[\phi] =  \left\{
\begin{array}{ll}
u_{m-k} & \mbox{when $1\leq k \leq m$}, \\
0   &  \mbox{when $k > m$},
\end{array}
\right. \\
\mathscr{D}_m [\phi \psi] &=& \sum_{k=0}^{m} \mathscr{D}_k[\phi] \; \mathscr{D}_{m-k}[\psi] = \sum_{k=0}^{m}u_k \; w_{m-k} \nonumber\\
&=& \sum_{k=0}^{m} \mathscr{D}_k[\psi] \; \mathscr{D}_{m-k}[\phi] = \sum_{k=0}^{m}w_k \; u_{m-k},\\
\mathscr{D}_m[\phi^{n+1}] &=& \sum_{k=0}^{m} \mathscr{D}_k[\phi] \; \mathscr{D}_{m-k}[\phi^n] = \sum_{k=0}^{m} u_k \;  \mathscr{D}_{m-k}[\phi^n].
\end{eqnarray}
Besides,  it holds
\begin{equation}
\mathscr{D}_m[\phi^n] =\sum_{r_1 =0}^{m} u_{m-r_1}\sum_{r_2=0}^{r_1}u_{r_1-r_2}\sum_{r_3=0}^{r_2}u_{r_2-r_3}\cdots\sum_{r_{n-1}=0}^{r_{n-2}}u_{r_{n-2}-r_{n-1}} u_{r_{n-1}},
\end{equation}
and
\begin{equation}
\mathscr{D}_m[f(x)\; \phi + g(x)\; \psi] =  f(x) \mathscr{D}_m[\phi] + g(x) \mathscr{D}_m[\psi]
\end{equation}
for arbitrary function $f$ and $g$ independent of $q$,  and
\begin{equation}
\mathscr{D}_m\left[{\mathscr L}[\phi]\right] = {\mathscr L}[\mathscr{D}_m[\phi]] ={\mathscr L}[u_m]
\end{equation}
for a linear operator $\mathscr L$, respectively.  In addition, it holds the recursion formulas
\begin{eqnarray}
{\mathscr D}_0\left[e^{\alpha \phi}\right] &=& e^{\alpha u_0}, \\
{\mathscr D}_m \left[e^{\alpha \phi}\right] &=& \alpha \sum_{k=0}^{m-1}\left( 1-\frac{k}{m}\right) u_{m-k} \; {\mathscr D}_k\left[e^{\alpha \phi}\right]; \\
{\mathscr D}_0[\sin\phi] &=& \sin(u_0),  \;\;\;  {\mathscr D}_0[\cos\phi] = \cos(u_0),\\
 {\mathscr D}_m\left[ \sin\phi\right]
 &=&\sum_{k=0}^{m-1}\left( 1-\frac{k}{m}\right) u_{m-k}\; {\mathscr D}_k\left[ \cos \phi \right], \\
{\mathscr D}_m\left[ \cos\phi\right]
 &=& -\sum_{k=0}^{m-1}\left( 1-\frac{k}{m}\right) u_{m-k}\; {\mathscr D}_k\left[ \sin \phi \right]
\end{eqnarray}
for $m\geq 1$.    In general,  it holds the recursion formulas
\begin{eqnarray}
{\mathscr D}_0[f(\phi)] &=& f(u_0), \\
{\mathscr D}_m[f(\phi)] &=& \sum_{k=0}^{m-1}\left(1-\frac{k}{m}\right) u_{m-k} \; {\mathscr D}_k\left[ f'(\phi)\right] \nonumber\\
&=& \sum_{k=0}^{m-1}\left(1-\frac{k}{m}\right) u_{m-k} \; \frac{\partial \left\{{\mathscr D}_k[f(\phi)] \right\}}{\partial u_0}
\end{eqnarray}
for $m\geq 1$, where $f(\phi)$ is a smooth function.

For detailed derivation of these properties, please refer to Liao \cite{Liao2009}  and  \S 4.2 of Liao's book \cite{liaobook2}.

Using the above properties, one can derive some other formulas.  For example,  it holds
\begin{eqnarray}
{\mathscr D}_{m+1}\left[ \sin(q \phi)\right] &=&  {\mathscr D}_{m+1}\left[  q \; \left\{ \frac{\sin(q \phi)}{q} \right\}\right]  \nonumber\\
&=& {\mathscr D}_1[q] {\mathscr D}_m \left[  \frac{\sin(q \phi)}{q} \right] \nonumber\\
&=&{\mathscr D}_m \left[  \frac{\sin(q \phi)}{q} \right].
\end{eqnarray}
So, using the recursion formulas mentioned above, one can get ${\mathscr D}_m[q^{-1} \sin(q\phi)]$.




\bibliographystyle{elsarticle-num}
\bibliography{InverseOperators}

\end{document}